\documentclass[11pt]{amsart}
\usepackage{amsmath,amssymb,mathrsfs}
\newtheorem{theorem}{Theorem}
\newtheorem{proposition}[theorem]{Proposition}
\newtheorem{lemma}[theorem]{Lemma}
\newtheorem{corollary}[theorem]{Corollary}

\begin{document}

\title[Yamabe problem on manifolds with boundary]{An existence theorem for the Yamabe problem on manifolds with boundary}
\author{Simon Brendle and Szu-Yu Sophie Chen}
\address{Department of Mathematics \\ Stanford University \\ Stanford, CA 94305}
\address{School of Mathematics \\ Institute for Advanced Study \\ Princeton, NJ 08540}
\thanks{The first author was partially supported by a Brown foundation fellowship and by the National Science Foundation under grant DMS-0905628. The second author was supported in part by the Miller Institute for Basic Research in Science and the National Science Foundation.}
\maketitle

\section{Introduction}

The Yamabe problem, solved by Trudinger \cite{Trudinger}, Aubin \cite{Aubin}, and Schoen \cite{Schoen}, asserts that any Riemannian metric on a closed manifold is conformal to a metric with constant scalar curvature. Escobar \cite{Escobar1}, \cite{Escobar2} has studied analogous questions on manifolds with boundary. To fix notation, let $(M,g)$ be a compact Riemannian manifold of dimension $n \geq 3$ with boundary $\partial M$. We denote by $R_g$ the scalar curvature of $(M,g)$ and by $\kappa_g$ the mean curvature of the boundary $\partial M$. There are two natural ways to extend the Yamabe problem to manifolds with boundary:
\begin{itemize}
\item[(a)] Find a metric $\tilde{g}$ in the conformal class of $g$ such that $R_{\tilde{g}}$ is constant and $\kappa_{\tilde{g}} = 0$.
\item[(b)] Find a metric $\tilde{g}$ in the conformal class of $g$ such that $R_{\tilde{g}} = 0$ and $\kappa_{\tilde{g}}$ is constant.
\end{itemize}
The boundary value problem (a) was first proposed by Escobar \cite{Escobar1}. The boundary value problem (b) is studied in \cite{Escobar2} and \cite{Marques}.

In this paper, we focus on the boundary value problem (a). The solvability of (a) is equivalent to the existence of a critical point of the Yamabe functional. This functional is defined by 
\[E_g(u) = \frac{\int_M (\frac{4(n-1)}{n-2} \, |du|_g^2 + R_g \, u^2) \, d\text{\rm vol}_g + \int_{\partial M} 2\kappa_g \, u^2 \, d\sigma_g}{\big ( \int_M u^{\frac{2n}{n-2}} \, d\text{\rm vol}_g \big )^{\frac{n-2}{n}}},\] 
where $u$ is a smooth positive function on $M$. Moreover, the Yamabe constant is defined as 
\[Y(M,\partial M,g) = \inf_{0 < u \in C^\infty(M)} E_g(u).\] 
It is well known that $Y(M,\partial M,g)$ is invariant under a conformal change of the metric $g$. Moreover, $Y(M,\partial M,g) \leq Y(S_+^n,\partial S_+^n)$, where $Y(S_+^n,\partial S_+^n)$ denotes the Yamabe constant of the hemisphere $S_+^n$ equipped with the standard metric.

Proving the existence of a minimizer for the functional $E_g$ is a difficult problem, as $E_g$ does not satisfy the Palais-Smale condition. The following existence result was established by Escobar \cite{Escobar1}. 

\begin{theorem}[J.~Escobar \cite{Escobar1}]
\label{escobar.theorem}
If $Y(M,\partial M,g) < Y(S_+^n,\partial S_+^n)$, then there exists a metric $\tilde{g}$ in the conformal class of $g$ such that $R_{\tilde{g}}$ is constant and $\kappa_{\tilde{g}}$ is equal to $0$.
\end{theorem}

Theorem \ref{escobar.theorem} should be compared to Aubin's existence theorem for the Yamabe problem on manifolds without boundary (cf. \cite{Aubin}).

In dimension $3 \leq n \leq 5$, Escobar showed that $Y(M,\partial M,g) < Y(S_+^n,\partial S_+^n)$ unless $M$ is conformally equivalent to the hemisphere $S_+^n$. In dimension $n \geq 6$, Escobar was able to verify the inequality $Y(M,\partial M,g) < Y(S_+^n,\partial S_+^n)$ under the assumption that the boundary $\partial M$ is not umbilic. 

Therefore, it remains to consider the case that $n \geq 6$ and $\partial M$ is umbilic. For abbreviation, we put $d = [\frac{n-2}{2}]$. As in \cite{Brendle2}, we denote by $\mathcal{Z}$ the set of all points $p \in M$ such that 
\[\limsup_{x \to p} d(p,x)^{2-d} \, |W_g(x)| = 0,\] 
where $W_g$ denotes the Weyl tensor of $(M,g)$. In other words, a point $p \in M$ belongs to $\mathcal{Z}$ if and only if $D^m W_g(p) = 0$ for all $m \in \{0,1,\hdots,d-2\}$. Note that the set $\mathcal{Z}$ is invariant under a conformal change of the metric.

The following is the main result of this paper:

\begin{theorem}
\label{main.theorem}
Let $(M,g)$ be a compact Riemannian manifold of dimension $n \geq 6$ with umbilic boundary $\partial M$. Moreover, let $p \in \partial M$ be an arbitrary point on the boundary of $M$. If $p \notin \mathcal{Z}$, then $Y(M,\partial M,g) < Y(S_+^n,\partial S_+^n)$. Consequently, there exists a metric $\tilde{g}$ in the conformal class of $g$ such that $R_{\tilde{g}}$ is constant and $\kappa_{\tilde{g}}$ is equal to $0$.
\end{theorem}

n case $p \in \mathcal{Z}$, we are able to show that $Y(M,\partial M,g) < Y(S_+^n,\partial S_+^N)$, provided that a certain asymptotically flat manifold has positive ADM mass (see Theorem \ref{degenerate.case} below).

We now give an outline of the proof of Theorem \ref{main.theorem}. By a theorem of Marques \cite{Marques}, we may work in conformal Fermi coordinates around $p$. We define a function $u_\varepsilon$ by 
\begin{equation} 
\label{definition.of.u}
u_\varepsilon(x) = \Big ( \frac{\varepsilon}{\varepsilon^2 + |x|^2} \Big )^{\frac{n-2}{2}}. 
\end{equation}
The function $u_\varepsilon$ satisfies 
\[\Delta u_\varepsilon = -n(n-2) \, u_\varepsilon^{\frac{n+2}{n-2}}\] 
and 
\[u_\varepsilon \, \partial_i \partial_k u_\varepsilon - \frac{n}{n-2} \, \partial_i u_\varepsilon \, \partial_k u_\varepsilon = \frac{1}{n} \, \Big ( u_\varepsilon \, \Delta u_\varepsilon - \frac{n}{n-2} \, |du_\varepsilon|^2 \Big ) \, \delta_{ik}.\] 
These identities reflect the fact that the metric $u_\varepsilon^{\frac{4}{n-2}} \, \delta_{ik}$ is Einstein.

We then consider a sum of the form $u_\varepsilon + w$, where $u_\varepsilon$ is given by (\ref{definition.of.u}) and $w$ is a correction term. This function is only defined in a small neighborhood of the point $p$. In order to extend the testfunction to all of $M$, we glue the function $u_\varepsilon + w$ to the Greens function of the conformal Laplacian with pole at $p$.

In order to show that the resuling testfunction has Yamabe energy less than $Y(S_+^n,\partial S_+^n)$, we make extensive use of techniques developed in \cite{Brendle2} (see also \cite{Brendle3}, \cite{Brendle-Marques}, \cite{Brendle-survey}). In \cite{Brendle2}, these techniques were used to prove a convergence theorem for the parabolic Yamabe flow in dimension $n \geq 6$. The convergence of the Yamabe flow in dimension $3 \leq n \leq 5$ was shown in \cite{Brendle1}.

\section{Auxiliary results}

In this section, we consider the halfspace $\mathbb{R}_+^n = \{x \in \mathbb{R}^n: x_n \geq 0\}$. 
Moreover, we assume that $H_{ik}(x)$ is a trace-free symmetric two-tensor on $\mathbb{R}_+^n$ which satisfies the following conditions: 
\begin{itemize}
\item At each point $x \in \mathbb{R}_+^n$, we have $H_{in}(x) = 0$ for all $i \in \{1,\hdots,n\}$.
\item At each point $x \in \partial \mathbb{R}_+^n$, we have $\sum_{k=1}^n H_{ik}(x) \, x_k = 0$ for all $i \in \{1,\hdots,n\}$.
\item At each point $x \in \partial \mathbb{R}_+^n$, we have $\partial_n H_{ik}(x) = 0$ for all $i,k \in \{1,\hdots,n\}$.
\end{itemize}
Finally, we assume that the components $H_{ik}(x)$ are polynomials of the form 
\[H_{ik}(x) = \sum_{2 \leq |\alpha| \leq d} h_{ik,\alpha} \, x^\alpha,\] 
where the sum is taken over all multi-indices $\alpha$ of length $2 \leq |\alpha| \leq d$.

As in \cite{Brendle2}, we define 
\[A_{ik} = \sum_{m=1}^n \partial_i \partial_m H_{mk} + \sum_{m=1}^n \partial_m \partial_k H_{im} - \Delta H_{ik} - \frac{1}{n-1} \sum_{m,p=1}^n \partial_m \partial_p H_{mp} \, \delta_{ik}\]
and
\begin{align*}
Z_{ijkl}
&= \partial_i \partial_k H_{jl} - \partial_i \partial_l H_{jk} - \partial_j \partial_k H_{il} + \partial_j \partial_l H_{ik} \\
&+ \frac{1}{n-2} \, (A_{jl} \, \delta_{ik} - A_{jk} \, \delta_{il} - A_{il} \, \delta_{jk} + A_{ik} \, \delta_{jl}).
\end{align*} 
Note that 
\[\partial_l Z_{ijkl} = \frac{n-3}{n-2} \, (\partial_i A_{jk} - \partial_j A_{ik}).\]

\begin{lemma}
\label{Weyl.tensor.1}
We have $A_{in}(x) = 0$ for all points $x \in \partial \mathbb{R}_+^n$ and all indices $i \in \{1,\hdots,n-1\}$.
\end{lemma}

\textbf{Proof.}
Note that $\partial_n H_{im}(x) = 0$ for all points $x \in \partial \mathbb{R}_+^n$ and all $i,m \in \{1,\hdots,n-1\}$. If we differentiate this identity in tangential direction, we obtain 
\[A_{in}(x) = \sum_{m=1}^{n-1} \partial_m \partial_n H_{im}(x) = 0\] 
for all points $x \in \partial \mathbb{R}_+^n$ and all $i \in \{1,\hdots,n-1\}$. \\

\begin{lemma}
\label{Weyl.tensor.2}
Assume that $Z_{ijkl}(x) = 0$ for all points $x \in \partial \mathbb{R}_+^n$ and all $i,j,k,l \in \{1,\hdots,n\}$. Then $H_{ik}(x) = A_{ik}(x) = 0$ for all $x \in \partial \mathbb{R}_+^n$ and all $i,k \in \{1,\hdots,n-1\}$.
\end{lemma}

\textbf{Proof.} 
We define 
\begin{align*} 
\hat{A}_{ik}(x) 
&= \sum_{m=1}^{n-1} \partial_i \partial_m H_{mk}(x) + \sum_{m=1}^{n-1} \partial_m \partial_k H_{im}(x) \\ 
&- \sum_{m=1}^{n-1} \partial_m \partial_m H_{ik}(x) - \frac{1}{n-2} \sum_{m,p=1}^n \partial_m \partial_p H_{mp}(x) \, \delta_{ik} 
\end{align*}
for all points $x \in \partial \mathbb{R}_+^n$ and all indices $i,k \in \{1,\hdots,n-1\}$. By assumption, we have 
\[\partial_n \partial_n H_{ik}(x) + \frac{1}{n-2} \, (A_{ik}(x) + A_{nn}(x) \, \delta_{ik}) = Z_{inkn}(x) = 0\] 
for all points $x \in \partial \mathbb{R}_+^n$ and all indices $i,k \in \{1,\hdots,n-1\}$. This implies  
\begin{align*} 
\hat{A}_{ik}(x) 
&= A_{ik}(x) + \partial_n \partial_n H_{ik}(x) - \frac{1}{(n-1)(n-2)} \sum_{m,p=1}^{n-1} \partial_m \partial_p H_{mp}(x) \, \delta_{ik} \\ 
&= A_{ik}(x) + \partial_n \partial_n H_{ik}(x) + \frac{1}{n-2} \, A_{nn}(x) \, \delta_{ik} \\ 
&= \frac{n-3}{n-2} \, A_{ik}(x) 
\end{align*} 
for all points $x \in \partial \mathbb{R}_+^n$ and all indices $i,k \in \{1,\hdots,n-1\}$. Hence, we obtain 
\begin{align*}
&\partial_i \partial_k H_{jl}(x) - \partial_i \partial_l H_{jk}(x) - \partial_j \partial_k H_{il}(x) + \partial_j \partial_l H_{ik}(x) \\ 
&= -\frac{1}{n-2} \, (A_{jl}(x) \, \delta_{ik} - A_{jk}(x) \, \delta_{il} - A_{il}(x) \, \delta_{jk} + A_{ik}(x) \, \delta_{jl}) \\ 
&= -\frac{1}{n-3} \, (\hat{A}_{jl}(x) \, \delta_{ik} - \hat{A}_{jk}(x) \, \delta_{il} - \hat{A}_{il}(x) \, \delta_{jk} + \hat{A}_{ik}(x) \, \delta_{jl}) 
\end{align*} 
for all points $x \in \partial \mathbb{R}_+^n$ and all indices $i,j,k,l \in \{1,\hdots,n-1\}$. Using Proposition 7 in \cite{Brendle2}, we conclude that $H_{ik}(x) = 0$ for all $x \in \partial \mathbb{R}_+^n$ and all $i,k \in \{1,\hdots,n-1\}$. This implies $A_{ik}(x) = \frac{n-2}{n-3} \, \hat{A}_{ik}(x) = 0$ for all $x \in \partial \mathbb{R}_+^n$ and all $i,k \in \{1,\hdots,n-1\}$. \\

\begin{proposition}
\label{Weyl.tensor.3}
Assume that $Z_{ijkl}(x) = 0$ for all $x \in \mathbb{R}_+^n$ and all $i,j,k,l \in \{1,\hdots,n\}$. Then $H_{ik}(x) = 0$ for all $x \in \mathbb{R}_+^n$ and all $i,k \in \{1,\hdots,n-1\}$.
\end{proposition}

\textbf{Proof.} 
Without loss of generality, we may assume that $H_{ik}(x)$ is homogenous of degree $d' \geq 2$. By assumption, we have $Z_{ijkl}(x) = 0$ for all $x \in \mathbb{R}_+^n$ and all $i,j,k,l \in \{1,\hdots,n\}$. This implies 
\[\partial_i A_{jk}(x) - \partial_j A_{ik}(x) = \frac{n-2}{n-3} \sum_{l=1}^n \partial_l Z_{ijkl}(x) = 0\] 
for all $x \in \mathbb{R}_+^n$ and all $i,j,k \in \{1,\hdots,n\}$. We next define 
\[\varphi(x) = \frac{1}{d'(d'-1)} \sum_{i,k=1}^n A_{ik}(x) \, x_i \, x_k\] 
for all $x \in \mathbb{R}_+^n$. Clearly, $\varphi(x)$ is a homogeneous polynomial of degree $d'$. Moreover, we have 
\[\partial_k \varphi(x) = \frac{1}{d'-1} \sum_{i=1}^n A_{ik}(x) \, x_i\] 
for all $x \in \mathbb{R}_+^n$ and all $k \in \{1,\hdots,n\}$. This implies 
\[\partial_i \partial_k \varphi(x) = A_{ik}(x)\] 
for all $x \in \mathbb{R}_+^n$ and all $i,k \in \{1,\hdots,n\}$. Using the identity $\sum_{i=1}^{n-1} H_{ii}(x) = 0$,we obtain 
\[\partial_n \partial_n \varphi(x) = A_{nn}(x) = -\frac{1}{n-1} \sum_{i=1}^{n-1} A_{ii}(x) = -\frac{1}{n-1} \sum_{i=1}^{n-1} \partial_i \partial_i \varphi(x)\] 
for all $x \in \mathbb{R}_+^n$. By Lemma \ref{Weyl.tensor.1}, we have $\partial_n \varphi(x) = \frac{1}{d'-1} \sum_{i=1}^{n-1} A_{in}(x) \, x_i = 0$ for all $x \in \partial \mathbb{R}_+^n$. Moreover, it follows from Lemma \ref{Weyl.tensor.2} that $\varphi(x) = \frac{1}{d'(d'-1)} \sum_{i,k=1}^{n-1} A_{ik}(x) \, x_i \, x_k = 0$ for all $x \in \partial \mathbb{R}_+^n$. Putting these facts together, we conclude that $\varphi(x) = 0$ for all $x \in \mathbb{R}_+^n$. 

Therefore, we have $A_{ik}(x) = \partial_i \partial_k \varphi(x) = 0$ for all $x \in \mathbb{R}_+^n$ and all $i,k \in \{1,\hdots,n\}$. This implies 
\[\partial_n \partial_n H_{ik}(x) = Z_{inkn}(x) - \frac{1}{n-2} \, (A_{ik}(x) + A_{nn}(x) \, \delta_{ik}) = 0\] 
for all $x \in \mathbb{R}_+^n$ and all $i,k \in \{1,\hdots,n-1\}$. On the other hand, we know that $H_{ik}(x) = \partial_n H_{ik}(x) = 0$ for all $x \in \partial \mathbb{R}_+^n$ and all $i,k \in \{1,\hdots,n-1\}$. From this, it follows that $H_{ik}(x) = 0$ for all $x \in \mathbb{R}_+^n$ and all $i,k \in \{1,\hdots,n-1\}$. This completes the proof of Proposition \ref{Weyl.tensor.3}. \\

For each $r > 0$, we denote by $U_r \subset \mathbb{R}_+^n$ the open ball of radius $\frac{r}{4}$ centered at the point $(0,\hdots,0,\frac{3r}{2})$. Moreover, let $u_\varepsilon: \mathbb{R}_+^n \to \mathbb{R}$ be defined by (\ref{definition.of.u}).

\begin{proposition}
\label{metric.bounded.by.Weyl.tensor}
There exists a constant $K_1$, depending only on $n$, such that
\[\sum_{2 \leq |\alpha| \leq d} \sum_{i,k=1}^n |h_{ik,\alpha}|^2 \, r^{2|\alpha|-4+n} \leq K_1 \, \int_{U_r} \sum_{i,j,k,l=1}^n |Z_{ijkl}(x)|^2 \, dx\] 
for all $r > 0$.
\end{proposition}

\textbf{Proof.}
It follows from Proposition \ref{Weyl.tensor.3} that the assertion holds for $r = 1$. The general case follows by scaling. \\

\begin{proposition}
\label{metric.bounded.by.E}
Let $V$ be a smooth vector field on $\mathbb{R}_+^n$. Moreover, let 
\[T_{ik} = H_{ik} - \partial_i V_k - \partial_k V_i + \frac{2}{n} \, \text{\rm div} \, V \, \delta_{ik}\]
and
\begin{align*}
Q_{ik,l}
&= u_\varepsilon \, \partial_l T_{ik} - \frac{2}{n-2} \, \partial_i u_\varepsilon \, T_{kl} - \frac{2}{n-2} \, \partial_k u_\varepsilon \, T_{il} \\
&+ \frac{2}{n-2} \sum_{p=1}^n \partial_p u_\varepsilon \, T_{ip} \, \delta_{kl} + \frac{2}{n-2} \sum_{p=1}^n \partial_p u_\varepsilon \, T_{kp} \, \delta_{il}. 
\end{align*}
Then there exists a constant $K_2$, depending only on $n$, such that 
\[\sum_{2 \leq |\alpha| \leq d} \sum_{i,k=1}^n |h_{ik,\alpha}|^2 \, \varepsilon^{n-2} \, r^{2|\alpha|+2-n} \leq K_2 \int_{(B_{2r}(0) \setminus B_r(0)) \cap \mathbb{R}_+^n} \sum_{i,k,l=1}^n |Q_{ik,l}(x)|^2 \, dx\] 
for all $r \geq \varepsilon$.
\end{proposition}

\textbf{Proof.}
In \cite{Brendle2}, the first author showed that 
\begin{align*} 
\frac{1}{4} \, \sum_{i,j,k,l=1}^n |Z_{ijkl}|^2 
&= \sum_{i,j,k,l=1}^n \partial_j(u_\varepsilon^{-1} \, Q_{ik,l}) \, Z_{ijkl} \\ 
&+ \frac{2}{n-2} \, \sum_{i,j,k,l=1}^n u_\varepsilon^{-2} \, \partial_k u_\varepsilon \, Q_{il,j} \, Z_{ijkl} 
\end{align*} 
(cf. \cite{Brendle2}, p.~555). Let us fix a smooth cut-off function $\eta: \mathbb{R}^n \to [0,1]$ such that $\eta(x) = 1$ for $x \in U_1$ and $\eta(x) = 0$ for $x \notin (B_2(0) \setminus B_1(0)) \cap \mathbb{R}_+^n$. In particular, we have $\eta(x) = 0$ for all $x \in \partial \mathbb{R}_+^n$. Integration by parts gives
\begin{align*}
&\int_{\mathbb{R}_+^n} \frac{1}{4} \, \sum_{i,j,k,l=1}^n |Z_{ijkl}(x)|^2 \, \eta(x/r) \, dx \\ 
&= -\int_{\mathbb{R}_+^n} \sum_{i,j,k,l=1}^n u_\varepsilon(x)^{-1} \, Q_{ik,l}(x) \, \partial_j \big [ Z_{ijkl}(x) \, \eta(x/r) \big ] \, dx \\ 
&+ \int_{\mathbb{R}_+^n} \frac{2}{n-2} \, \sum_{i,j,k,l=1}^n u_\varepsilon(x)^{-2} \, \partial_k u_\varepsilon(x) \, Q_{il,j}(x) \, Z_{ijkl}(x) \, \eta(x/r) \, dx. 
\end{align*}
Using H\"older's inequality, we obtain 
\begin{align*}
&\int_{U_r} \sum_{i,j,k,l=1}^n |Z_{ijkl}(x)|^2 \, dx \\ 
&\leq K_3 \, \varepsilon^{-\frac{n-2}{2}} \, r^{n-3} \, \bigg ( \sum_{2 \leq |\alpha| \leq d} \sum_{i,k=1}^n |h_{ik,\alpha}|^2 \, r^{2|\alpha|-4+n} \bigg )^{\frac{1}{2}} \\ 
&\hspace{10mm} \cdot \bigg ( \int_{(B_{2r} \setminus B_r(0)) \cap \mathbb{R}_+^n} \sum_{i,k,l=1}^n |Q_{ik,l}(x)|^2 \, dx \bigg )^{\frac{1}{2}} 
\end{align*} 
for all $r \geq \varepsilon$. Here, $K_3$ is a positive constant that depends only on $n$. On the other hand, it follows from Proposition \ref{metric.bounded.by.Weyl.tensor} that 
\[\sum_{2 \leq |\alpha| \leq d} \sum_{i,k=1}^n |h_{ik,\alpha}|^2 \, r^{2|\alpha|-4+n} \leq K_1 \int_{U_r} \sum_{i,j,k,l=1}^n |Z_{ijkl}(x)|^2 \, dx.\]
Putting these facts together, the assertion follows. \\

\begin{corollary}
\label{good.term}
Let $V$ be a smooth vector field on $\mathbb{R}_+^n$. Moreover, let
\[T_{ik} = H_{ik} - \partial_i V_k - \partial_k V_i + \frac{2}{n} \, \text{\rm div} \, V \, \delta_{ik}\]
and
\begin{align*}
Q_{ik,l}
&= u_\varepsilon \, \partial_l T_{ik} - \frac{2}{n-2} \, \partial_i u_\varepsilon \, T_{kl} - \frac{2}{n-2} \, \partial_k u_\varepsilon \, T_{il} \\
&+ \frac{2}{n-2} \sum_{p=1}^n \partial_p u_\varepsilon \, T_{ip} \, \delta_{kl} + \frac{2}{n-2} \sum_{p=1}^n \partial_p u_\varepsilon \, T_{kp} \, \delta_{il}.
\end{align*}
Then there exists a constant $K_4$, depending only on $n$, such that
\begin{align*}
&\sum_{2 \leq |\alpha| \leq d} \sum_{i,k=1}^n |h_{ik,\alpha}|^2 \, \varepsilon^{n-2} \, \int_{B_\delta(0) \cap \mathbb{R}_+^n} (\varepsilon + |x|)^{2|\alpha|+2-2n} \, dx \\ 
&\leq K_4 \int_{B_\delta(0) \cap \mathbb{R}_+^n} \sum_{i,k,l=1}^n |Q_{ik,l}(x)|^2 \, dx 
\end{align*}
for all $\delta \geq 2\varepsilon$.
\end{corollary}

\section{The main estimate}

We now describe the construction of the test function. Let $(M,g)$ be a compact Riemannian manifold of dimension $n \geq 6$ with umbilic boundary $\partial M$. After changing the metric conformally, we may assume that $\partial M$ is totally geodesic.

Let us fix a point $p \in \partial M$, and let $(x_1,\hdots,x_n)$ denote the Fermi coordinates around $p$. In these coordinates, the metric has the following properties: 
\begin{itemize}
\item At each point $x \in \mathbb{R}_+^n$, we have $g_{in}(x) = \delta_{in}$ for all $i \in \{1,\hdots,n\}$.
\item At each point $x \in \partial \mathbb{R}_+^n$, we have $\sum_{k=1}^n g_{ik}(x) \, x_k = x_i$ for all $i \in \{1,\hdots,n\}$.
\item At each point $x \in \partial \mathbb{R}_+^n$, we have $\partial_n g_{ik}(x) = 0$ for all $i,k \in \{1,\hdots,n\}$.
\end{itemize}
By a theorem of Marques, there exists a system of conformal Fermi coordinates around $p$ (see \cite{Marques}, Proposition 3.1). Hence, after performing a conformal change of the metric, we may assume that $\det g(x) = 1 + O(|x|^{2d+2})$, where $d = [\frac{n-2}{2}]$. 

In the next step, we write $g(x) = \exp(h(x))$, where $h(x)$ is a smooth function taking values in the space of symmetric $n \times n$ matrices. This function has the following properties: 
\begin{itemize}
\item At each point $x \in \mathbb{R}_+^n$, we have $h_{in}(x) = 0$ for all $i \in \{1,\hdots,n\}$.
\item At each point $x \in \partial \mathbb{R}_+^n$, we have $\sum_{k=1}^n h_{ik}(x) \, x_k = 0$ for all $i \in \{1,\hdots,n\}$.
\item At each point $x \in \partial \mathbb{R}_+^n$, we have $\partial_n h_{ik}(x) = 0$ for all $i,k \in \{1,\hdots,n\}$.
\end{itemize}
Moreover, we have $\text{\rm tr} \, h(x) = O(|x|^{2d+2})$. For abbreviation, we denote by 
\[H_{ik}(x) = \sum_{2 \leq |\alpha| \leq d} h_{ik,\alpha} \, x^\alpha\] 
the Taylor polynomial of order $d$ associated with the function $h_{ik}(x)$. Clearly, $H_{ik}(x)$ is a trace-free symmetric two-tensor on $\mathbb{R}_+^n$. Moreover, we have $h_{ik}(x) = H_{ik}(x) + O(|x|^{d+1})$.

Let us fix a non-negative smooth function such that $\chi(t) = 1$ for $t \leq \frac{4}{3}$ and $\chi(t) = 0$ for $t \geq \frac{5}{3}$. Given any $\delta > 0$, we define a cut-off function $\chi_\delta: \mathbb{R}^n \to \mathbb{R}$ by $\chi_\delta(x) = \chi(|x| / \delta)$. By Theorem \ref{existence.of.smooth.solution}, there exists a smooth vector field $V$ on $\mathbb{R}_+^n$ with the following properties:
\begin{itemize}
\item At each point $x \in \mathbb{R}_+^n$, we have 
\[\sum_{k=1}^n \partial_k \Big [ u_\varepsilon^{\frac{2n}{n-2}} \, (\chi_\delta \, H_{ik} - \partial_i V_k - \partial_k V_i + \frac{2}{n} \, \text{\rm div} \, V \, \delta_{ik}) \Big ] = 0\]
for all $i \in \{1, \hdots, n\}$. 
\item At each point $x \in \partial \mathbb{R}_+^n$, we have $V_n(x) = \partial_n V_i(x) = 0$ for all $i \in \{1,\hdots,n-1\}$.
\end{itemize}
By Corollary \ref{polynomial}, the vector field $V$ satisfies the estimate
\begin{equation} 
\label{estimate.for.V} 
|\partial^\beta V^{(\varepsilon,\delta)}(x)| \leq C \sum_{2 \leq |\alpha| \leq d} \sum_{i,k=1}^n |h_{ik,\alpha}| \, (\varepsilon + |x|)^{|\alpha|+1-|\beta|} 
\end{equation}
for every multi-index $\beta$ and all $x \in \mathbb{R}_+^n$. Here, $C$ is a positive constant that depends only on $n$ and $|\beta|$.

For abbreviation, we define 
\[S_{ik} = \partial_i V_k + \partial_k V_i - \frac{2}{n} \, \text{\rm div} \, V \, \delta_{ik},\] 
\[T_{ik} = H_{ik} - S_{ik},\] 
\begin{align*}
Q_{ik,l}
&= u_\varepsilon \, \partial_l T_{ik} - \frac{2}{n-2} \, \partial_i u_\varepsilon \, T_{kl} - \frac{2}{n-2} \, \partial_k u_\varepsilon \, T_{il} \\
&+ \frac{2}{n-2} \sum_{p=1}^n \partial_p u_\varepsilon \, T_{ip} \, \delta_{kl} + \frac{2}{n-2} \sum_{p=1}^n \partial_p u_\varepsilon \, T_{kp} \, \delta_{il}, 
\end{align*}
and 
\[w = \sum_{l=1} \partial_l u_\varepsilon \, V_l + \frac{n-2}{2n} \, u_\varepsilon \, \text{\rm div} \, V.\] 
By definition of $V$, we have 
\begin{equation} 
\label{elliptic.system.1}
\sum_{k=1}^n \partial_k(u_\varepsilon^{\frac{2n}{n-2}} \, T_{ik}) = 0 
\end{equation}
for all points $x \in B_\delta(0) \cap \mathbb{R}_+^n$ and all $i \in \{1, \hdots, n\}$. This implies 
\begin{equation} 
\label{elliptic.system.2}
\sum_{k=1}^n \Big ( u_\varepsilon \, \partial_k T_{ik} + \frac{2n}{n-2} \, \partial_k u_\varepsilon \, T_{ik} \Big ) = 0 
\end{equation}
for all points $x \in B_\delta(0) \cap \mathbb{R}_+^n$ and all $i \in \{1, \hdots, n\}$. The following result was established in \cite{Brendle2}:

\begin{proposition}[S.~Brendle \cite{Brendle2}]
\label{integration.by.parts.identity}
There exists a smooth vector field $\xi$ on $\mathbb{R}_+^n$ such that 
\begin{align*} 
&\frac{1}{4} \, u_\varepsilon^2 \sum_{i,k,l=1}^n \partial_l H_{ik} \, \partial_l H_{ik} - \frac{1}{2} \, u_\varepsilon^2 \sum_{i,k,l=1}^n \partial_k H_{ik} \, \partial_l H_{il} \\ 
&- 2 \, u_\varepsilon \sum_{i,k,l=1}^n \partial_k u_\varepsilon \, H_{ik} \, \partial_l H_{il} - \frac{2(n-1)}{n-2} \sum_{i,k,l=1}^n \partial_k u_\varepsilon \, \partial_l u_\varepsilon \, H_{ik} \, H_{il} \\ 
&- 2 \, u_\varepsilon \, w \sum_{i,k=1}^n \partial_i \partial_k H_{ik} + \frac{8(n-1)}{n-2} \sum_{i,k=1}^n \partial_i u_\varepsilon \, \partial_k w \, H_{ik} \notag \\ 
&- \frac{4(n-1)}{n-2} \, |dw|^2 + \frac{4(n-1)}{n-2} \, n(n+2) \, u_\varepsilon^{\frac{4}{n-2}} \, w^2 \\ 
&= \frac{1}{4} \sum_{i,k,l=1}^n Q_{ik,l} \, Q_{ik,l} + 2 \, u_\varepsilon^{\frac{2n}{n-2}} \, \sum_{i,k=1}^n T_{ik} \, T_{ik} + \text{\rm div} \, \xi 
\end{align*}
for all points $x \in B_\delta(0) \cap \mathbb{R}_+^n$.
\end{proposition}

The vector field $\xi$ can be expressed in terms of the tensor $H_{ik}$ and the vector field $V$ (cf. \cite{Brendle2}, Section 2). In the next step, we show that $\xi$ is tangential along $\partial \mathbb{R}_+^n$. To that end, we need the following lemma:

\begin{lemma}
\label{relations.along.boundary}
At each point $x \in B_\delta(0) \cap \partial \mathbb{R}_+^n$, we have 
\[S_{in}(x) = T_{in}(x) = 0\] 
and 
\[\partial_n S_{ik}(x) = \partial_n T_{ik}(x) = 0\] 
for all $i,k \in \{1, \hdots, n-1\}$. Moreover, we have $\partial_n S_{nn}(x) = \partial_n T_{nn}(x) = 0$ and $\partial_n w(x) = 0$.
\end{lemma}

\textbf{Proof.} 
By assumption, we have $V_n(x) = \partial_n V_i(x) = 0$ for all $x \in \partial \mathbb{R}_+^n$. 
This implies $S_{in}(x) = T_{in}(x) = 0$ for all $i \in \{1,\hdots,n-1\}$. This implies 
\[\sum_{k=1}^{n-1} \Big ( u_\varepsilon(x) \, \partial_k T_{kn}(x) + \frac{2n}{n-2} \, \partial_k u_\varepsilon(x) \, T_{kn}(x) \Big ) = 0.\] 
Using (\ref{elliptic.system.2}), we obtain 
\[u_\varepsilon(x) \, \partial_n T_{nn}(x) + \frac{2n}{n-2} \, \partial_n u_\varepsilon(x) \, T_{nn}(x) = 0.\] 
This implies $\partial_n T_{nn}(x) = 0$, hence $\partial_n S_{nn}(x) = 0$. Consequently, we have $\partial_n \partial_n V(x) = 0$. From this, the assertion follows easily. \\

\begin{lemma}
\label{xi.tangential.along.boundary}
We have $\xi_n(x) = 0$ for all points $x \in B_\delta(0) \cap \partial \mathbb{R}_+^n$.
\end{lemma}

\textbf{Proof.} 
The vector field $\xi$ satisfies 
\begin{align*} 
\xi_n 
&= -2 \sum_{k=1}^n u_\varepsilon \, w \, \partial_k H_{nk} + 2 \sum_{i=1}^n \partial_i (u_\varepsilon \, w) \, H_{in} \\ 
&+ \frac{1}{2} \, u_\varepsilon^2 \sum_{i,k=1}^n \partial_n S_{ik} \, H_{ik} - u_\varepsilon^2 \sum_{i,l=1}^n \partial_l S_{il} \, H_{in} - 2 \, u_\varepsilon \sum_{i,l=1}^n \partial_l u_\varepsilon \, S_{il} \, H_{in} \\ 
&+ u_\varepsilon \, w \sum_{k=1}^n \partial_k S_{nk} - \sum_{i=1}^n \partial_i (u_\varepsilon \, w) \, S_{in} \\ 
&- \frac{1}{4} \, u_\varepsilon^2 \sum_{i,k=1}^n \partial_n S_{ik} \, S_{ik} + \frac{1}{2} \, u_\varepsilon^2 \sum_{i,l=1}^n \partial_l S_{il} \, S_{in} + u_\varepsilon \sum_{i,l=1}^n \partial_l u_\varepsilon \, S_{il} \, S_{in} \\ 
&+ \frac{4(n-1)}{n-2} \sum_{i=1}^n \partial_i u_\varepsilon \, w \, S_{in} - \frac{4(n-1)}{n-2} \, w \, \partial_n w \\ 
&+ \frac{2}{n-2} \, u_\varepsilon \sum_{i,k=1}^n \partial_k u_\varepsilon \, T_{ik} \, T_{in} 
\end{align*} 
(see \cite{Brendle2}, Section 2). Using Lemma \ref{relations.along.boundary}, we conclude that $\xi_n(x) = 0$ for all $x \in B_\delta(0) \cap \partial \mathbb{R}_+^n$. \\

\begin{proposition}
\label{integrated.version}
We have 
\begin{align*} 
&\int_{B_\delta(0) \cap \mathbb{R}_+^n} \bigg [ \frac{1}{4} \, u_\varepsilon^2 \sum_{i,k,l=1}^n \partial_l H_{ik} \, \partial_l H_{ik} - \frac{1}{2} \, u_\varepsilon^2 \sum_{i,k,l=1}^n \partial_k H_{ik} \, \partial_l H_{il} \bigg ] \\ 
&- \int_{B_\delta(0) \cap \mathbb{R}_+^n} \bigg [ 2 \, u_\varepsilon \sum_{i,k,l=1}^n \partial_k u_\varepsilon \, H_{ik} \, \partial_l H_{il} + \frac{2(n-1)}{n-2} \sum_{i,k,l=1}^n \partial_k u_\varepsilon \, \partial_l u_\varepsilon \, H_{ik} \, H_{il} \bigg ] \\ 
&- \int_{B_\delta(0) \cap \mathbb{R}_+^n} \bigg [ 2 \, u_\varepsilon \, w \sum_{i,k=1}^n \partial_i \partial_k H_{ik} - \frac{8(n-1)}{n-2} \sum_{i,k=1}^n \partial_i u_\varepsilon \, \partial_k w \, H_{ik} \bigg ] \\ 
&- \int_{B_\delta(0) \cap \mathbb{R}_+^n} \bigg [ \frac{4(n-1)}{n-2} \, |dw|^2 - \frac{4(n-1)}{n-2} \, n(n+2) \, u_\varepsilon^{\frac{4}{n-2}} \, w^2 \bigg ] \\ 
&\geq 2\lambda \sum_{2 \leq |\alpha| \leq d} \sum_{i,k=1}^n |h_{ik,\alpha}|^2 \, \varepsilon^{n-2} \, \int_{B_\delta(0) \cap \mathbb{R}_+^n} (\varepsilon + |x|)^{2|\alpha|+2-2n} \, dx \\ 
&- C \sum_{2 \leq |\alpha| \leq d} \sum_{i,k=1}^n |h_{ik,\alpha}|^2 \, \delta^{2|\alpha|+2-n} \, \varepsilon^{n-2}. 
\end{align*}
for $\delta \geq 2\varepsilon$. Here, $\lambda$ and $C$ are positive constants that depend only on $n$.
\end{proposition}

\textbf{Proof.}
We consider the identity in Proposition \ref{integration.by.parts.identity} and integrate over $B_\delta(0) \cap \mathbb{R}_+^n$. By Corollary \ref{good.term}, we have 
\begin{align*}
&\int_{B_\delta(0) \cap \mathbb{R}_+^n} \sum_{i,k,l=1}^n Q_{ik,l} \, Q_{ik,l} \\ 
&\geq 8\lambda \sum_{2 \leq |\alpha| \leq d} \sum_{i,k=1}^n |h_{ik,\alpha}|^2 \, \varepsilon^{n-2} \, \int_{B_\delta(0) \cap \mathbb{R}_+^n} (\varepsilon + |x|)^{2|\alpha|+3-2n} \, dx, 
\end{align*}
where $\lambda = 1 / (8K_4)$ is a positive constant that depends only on $n$. Moreover, it follows from \ref{estimate.for.V} that 
\[|\xi(x)| \leq C \, \varepsilon^{n-2} \sum_{2 \leq |\alpha| \leq d} \sum_{i,k=1}^n |h_{ik,\alpha}| \, (\varepsilon + |x|)^{2|\alpha|+3-2n}\] 
for all $x \in \mathbb{R}_+^n$. Using Lemma \ref{xi.tangential.along.boundary} and the divergence theorem, we obtain 
\begin{align*} 
\int_{B_\delta(0) \cap \mathbb{R}_+^n} \text{\rm div} \, \xi 
&= \int_{\partial B_\delta(0) \cap \mathbb{R}_+^n} \sum_{i=1}^n \frac{x_i}{|x|} \, \xi_i - \int_{B_\delta(0) \cap \partial \mathbb{R}_+^n} \xi_n \\ 
&\leq C \, \varepsilon^{n-2} \sum_{2 \leq |\alpha| \leq d} \sum_{i,k=1}^n |h_{ik,\alpha}| \, \delta^{2|\alpha|+2-n}. 
\end{align*}
Putting these facts together, the assertion follows. \\

Finally, we need the following estimate for the scalar curvature $R_g$.

\begin{proposition}
\label{Taylor.expansion.of.scalar.curvature}
The scalar curvature $R_g$ satisfies the estimates
\begin{equation} 
\label{scalar.curvature.1}
\bigg | R_g - \sum_{i,k=1}^n \partial_i \partial_k H_{ik} \bigg | \leq C \sum_{2 \leq |\alpha| \leq d} \sum_{i,k=1}^n |h_{ik,\alpha}| \, |x|^{|\alpha|} + C \, |x|^{d-1} 
\end{equation} 
and 
\begin{align}
\label{scalar.curvature.2}
&\bigg | R_g - \sum_{i,k=1}^n \partial_i \partial_k h_{ik} + \sum_{i,k,l=1}^n \partial_k(H_{ik} \, \partial_l H_{il}) \notag \\ 
&\hspace{10mm} - \frac{1}{2} \sum_{i,k,l=1}^n \partial_k H_{ik} \, \partial_l H_{il} + \frac{1}{4} \sum_{i,k,l=1}^n \partial_l H_{ik} \, \partial_l H_{ik} \bigg | \notag \\ 
&\leq C \sum_{2 \leq |\alpha| \leq d} \sum_{i,k=1}^n |h_{ik,\alpha}|^2 \, |x|^{2|\alpha|} \\ 
&+ C \sum_{2 \leq |\alpha| \leq d} \sum_{i,k=1}^n |h_{ik,\alpha}| \, |x|^{|\alpha|+d-1} + C \, |x|^{2d} \notag
\end{align}
if $|x|$ is sufficiently small.
\end{proposition}

\textbf{Proof.} 
This follows easily from Proposition 25 in \cite{Brendle3} (see also Corollary 12 in \cite{Brendle2}, where geodesic normal coordinates are considered). \\

Our goal is to estimate the Yamabe energy of $u_\varepsilon + w$. To that end, we proceed in several steps:

\begin{proposition}
\label{estimate.for.test.function}
There exist positive constants $\lambda$, $C$, $\delta_0$ such that 
\begin{align*}
&\int_{B_\delta(0) \cap \mathbb{R}_+^n} \Big ( \frac{4(n-1)}{n-2} \, |d(u_\varepsilon + w)|_g^2 + R_g \, (u_\varepsilon + w)^2 \Big ) \\
&\leq \int_{B_\delta(0) \cap \mathbb{R}_+^n} \frac{4(n-1)}{n-2} \, |du_\varepsilon|^2 + \int_{B_\delta(0) \cap \mathbb{R}_+^n} \frac{4(n-1)}{n-2} \, n(n+2) \, u_\varepsilon^{\frac{4}{n-2}} \, w^2 \\ 
&+ \int_{\partial B_\delta(0) \cap \mathbb{R}_+^n} \sum_{i,k=1}^n \frac{x_i}{|x|} \, (u_\varepsilon^2 \, \partial_k h_{ik} - 2 \, u_\varepsilon \, \partial_k u_\varepsilon \, h_{ik}) \\ 
&- \lambda \sum_{2 \leq |\alpha| \leq d} \sum_{i,k=1}^n |h_{ik,\alpha}|^2 \, \varepsilon^{n-2} \, \int_{B_\delta(0) \cap \mathbb{R}_+^n} (\varepsilon + |x|)^{2|\alpha|+2-2n} \, dx \\ 
&+ C \sum_{2 \leq |\alpha| \leq d}  \sum_{i,k=1}^n |h_{ik,\alpha}| \, \delta^{|\alpha|+2-n} \, \varepsilon^{n-2} + C \, \delta^{2d+4-n} \, \varepsilon^{n-2} 
\end{align*}
if $0 < 2\varepsilon \leq \delta \leq \delta_0$. The constant $\lambda$ depends only on $n$. The constants $C$, $\delta_0$ depend on the background manifold $(M,g)$.
\end{proposition}

\textbf{Proof.}
Let us write 
\begin{align*}
&\frac{4(n-1)}{n-2} \, |d(u_\varepsilon + w)|_g^2 + R_g \, (u_\varepsilon + w)^2 \\
&= \frac{4(n-1)}{n-2} \, |du_\varepsilon|^2 + \frac{4(n-1)}{n-2} \, n(n+2) \, u_\varepsilon^{\frac{4}{n-2}} \, w^2 \\ 
&+ \frac{8(n-1)}{n-2} \, J^{(1)} + J^{(2)} + J^{(3)} + J^{(4)} + J^{(5)} + J^{(6)} + J^{(7)}, 
\end{align*} 
where 
\begin{align*} 
J^{(1)} &= \sum_{i=1}^n \partial_i u_\varepsilon \, \partial_i w \\ 
J^{(2)} &= -\frac{4(n-1)}{n-2} \sum_{i,k=1}^n \partial_i u_\varepsilon \, \partial_k u_\varepsilon \, h_{ik} + u_\varepsilon^2 \sum_{i,k=1}^n \partial_i \partial_k h_{ik} \\ 
J^{(3)} &= -u_\varepsilon^2 \sum_{i,k,l=1}^n \partial_k(H_{ik} \, \partial_l H_{il}) - 2 \, u_\varepsilon \sum_{i,k,l=1}^n \partial_k u_\varepsilon \, H_{ik} \, \partial_l H_{il} \\
J^{(4)} &= -\frac{1}{4} \sum_{i,k,l=1}^n u_\varepsilon^2 \, \partial_l H_{ik} \, \partial_l H_{ik} + \frac{1}{2} \sum_{i,k,l=1}^n u_\varepsilon^2 \, \partial_k H_{ik} \, \partial_l H_{il} \\ 
&+ 2 \, u_\varepsilon \sum_{i,k,l=1}^n \partial_k u_\varepsilon \, H_{ik} \, \partial_l H_{il} + \frac{2(n-1)}{n-2} \sum_{i,k,l=1}^n \partial_k u_\varepsilon \, \partial_l u_\varepsilon \, H_{ik} \, H_{il} \\ 
&+ 2 \, u_\varepsilon \, w \sum_{i,k=1}^n \partial_i \partial_k H_{ik} - \frac{8(n-1)}{n-2} \sum_{i,k=1}^n \partial_i u_\varepsilon \, \partial_k w \, H_{ik} \\ 
&+ \frac{4(n-1)}{n-2} \, |dw|^2 - \frac{4(n-1)}{n-2} \, n(n+2) \, u_\varepsilon^{\frac{4}{n-2}} \, w^2 \\ 
J^{(5)} &= \frac{4(n-1)}{n-2} \sum_{i,k=1}^n \bigg [ g^{ik} - \delta_{ik} + h_{ik} - \frac{1}{2} \, \sum_{l=1}^n H_{il} \, H_{kl} \bigg ] \, \partial_i u_\varepsilon \, \partial_k u_\varepsilon \\ 
&+ \bigg [ R_g - \sum_{i,k=1}^n \partial_i \partial_k h_{ik} + \sum_{i,k,l=1}^n \partial_k(H_{ik} \, \partial_l H_{il}) \\ 
&\hspace{10mm} - \frac{1}{2} \sum_{i,k,l=1}^n \partial_k H_{ik} \, \partial_l H_{il} + \frac{1}{4} \sum_{i,k,l=1}^n \partial_l H_{ik} \, \partial_l H_{ik} \bigg ] \, u_\varepsilon^2 
\end{align*} 
and 
\begin{align*}
J^{(6)} &= \frac{8(n-1)}{n-2} \sum_{i,k=1}^n (g^{ik} - \delta_{ik} + H_{ik}) \, \partial_i u_\varepsilon \, \partial_k w \\ 
&+ 2 \, \bigg [ R_g - \sum_{i,k=1}^n \partial_i \partial_k H_{ik} \bigg ] \, u_\varepsilon \, w \\ 
J^{(7)} &= R_g \, w^2 + \frac{4(n-1)}{n-2} \sum_{i,k=1}^n (g^{ik} - \delta_{ik}) \, \partial_i w \, \partial_k w. 
\end{align*}
It follows from the divergence theorem that 
\begin{align*} 
\int_{B_\delta(0) \cap \mathbb{R}_+^n} J^{(1)} 
&= \int_{B_\delta(0) \cap \mathbb{R}_+^n} \sum_{i=1}^n \partial_i \Big [ \partial_i u_\varepsilon \, w + \frac{(n-2)^2}{2} \, u_\varepsilon^{\frac{2n}{n-2}} \, V_i \Big ] \\ 
&= \int_{\partial B_\delta(0) \cap \mathbb{R}_+^n} \sum_{i=1}^n \frac{x_i}{|x|} \, \Big [ \partial_i u_\varepsilon \, w + \frac{(n-2)^2}{2} \, u_\varepsilon^{\frac{2n}{n-2}} \, V_i \Big ] \\ 
&\leq C \sum_{2 \leq |\alpha| \leq d} \sum_{i,k=1}^n |h_{ik,\alpha}| \, \delta^{|\alpha|+2-n} \, \varepsilon^{n-2}. 
\end{align*} 
We next observe that 
\begin{align*} 
&J^{(2)} - \sum_{i,k=1}^n \partial_i(u_\varepsilon^2 \, \partial_k h_{ik} - 2 \, u_\varepsilon \, \partial_k u_\varepsilon \, h_{ik}) \\ 
&= 2 \sum_{i,k=1}^n \Big ( u_\varepsilon \, \partial_i \partial_k u_\varepsilon - \frac{n}{n-2} \, \partial_i u_\varepsilon \, \partial_k u_\varepsilon \Big ) \, h_{ik} \\ 
&= \frac{2}{n} \, \Big ( u_\varepsilon \, \Delta u_\varepsilon - \frac{n}{n-2} \, |du_\varepsilon|^2 \Big ) \, \text{\rm tr} \, h \\ 
&\leq C \, \varepsilon^{n-2} \, (\varepsilon + |x|)^{2d+4-2n}. 
\end{align*} 
Using the divergence theorem, we obtain 
\begin{align*} 
&\int_{B_\delta(0) \cap \mathbb{R}_+^n} J^{(2)} \\ 
&\leq \sum_{i,k=1}^n \partial_i(u_\varepsilon^2 \, \partial_k h_{ik} - 2 \, u_\varepsilon \, \partial_k u_\varepsilon \, h_{ik}) + C \, \delta^{2d+4-n} \, \varepsilon^{n-2} \\ 
&\leq \int_{\partial B_\delta(0) \cap \mathbb{R}_+^n} \sum_{i,k=1}^n \frac{x_i}{|x|} \, (u_\varepsilon^2 \, \partial_k h_{ik} - 2 \, u_\varepsilon \, \partial_k u_\varepsilon \, h_{ik}) + C \, \delta^{2d+4-n} \, \varepsilon^{n-2}. 
\end{align*} 
Moreover, we have 
\begin{align*} 
\int_{B_\delta(0) \cap \mathbb{R}_+^n} J^{(3)} 
&= -\int_{B_\delta(0) \cap \mathbb{R}_+^n} \sum_{i,k,l=1}^n \partial_k(u_\varepsilon^2 \, H_{ik} \, \partial_l H_{il}) \\ 
&= -\int_{\partial B_\delta(0) \cap \mathbb{R}_+^n} \sum_{i,k,l=1}^n \frac{x_k}{|x|} \, u_\varepsilon^2 \, H_{ik} \, \partial_l H_{il} \\ 
&\leq C \sum_{2 \leq |\alpha| \leq d} \sum_{i,k=1}^n |h_{ik,\alpha}|^2 \, \delta^{2|\alpha|+2-n} \, \varepsilon^{n-2}. 
\end{align*} 
Using Proposition \ref{integrated.version}, we obtain 
\begin{align*} 
&\int_{B_\delta(0) \cap \mathbb{R}_+^n} J^{(4)} \\ 
&\leq -2\lambda \sum_{2 \leq |\alpha| \leq d} \sum_{i,k=1}^n |h_{ik,\alpha}|^2 \, \varepsilon^{n-2} \, \int_{B_\delta(0) \cap \mathbb{R}_+^n} (\varepsilon + |x|)^{2|\alpha|+2-2n} \, dx \\ 
&+ C \sum_{2 \leq |\alpha| \leq d} \sum_{i,k=1}^n |h_{ik,\alpha}|^2 \, \delta^{2|\alpha|+2-n} \, \varepsilon^{n-2}. 
\end{align*}
It remains to estimate the terms $J^{(5)}$, $J^{(6)}$, and $J^{(7)}$. Using Proposition \ref{Taylor.expansion.of.scalar.curvature}, we obtain the pointwise estimate 
\begin{align*} 
J^{(5)} + J^{(6)} + J^{(7)} 
&\leq C \sum_{2 \leq |\alpha| \leq d} \sum_{i,k=1}^n |h_{ik,\alpha}|^2 \, \varepsilon^{n-2} \, (\varepsilon + |x|)^{2|\alpha|+4-2n} \\ 
&+ C \sum_{2 \leq |\alpha| \leq d} \sum_{i,k=1}^n |h_{ik,\alpha}| \, \varepsilon^{n-2} \, (\varepsilon + |x|)^{|\alpha|+d+3-2n} \\ 
&+ C \, \varepsilon^{n-2} \, (\varepsilon + |x|)^{2d+4-2n} 
\end{align*}
for $x \in B_\delta(0) \cap \mathbb{R}_+^n$. Using Young's inequality, we deduce that 
\begin{align*} 
J^{(5)} + J^{(6)} + J^{(7)} 
&\leq \lambda \sum_{2 \leq |\alpha| \leq d} \sum_{i,k=1}^n |h_{ik,\alpha}|^2 \, \varepsilon^{n-2} \, (\varepsilon + |x|)^{2|\alpha|+2-2n} \\ 
&+ C \, \varepsilon^{n-2} \, (\varepsilon + |x|)^{2d+4-2n} 
\end{align*}
for $x \in B_\delta(0) \cap \mathbb{R}_+^n$. Integration over $B_\delta(0) \cap \mathbb{R}_+^n$ yields 
\begin{align*} 
&\int_{B_\delta(0) \cap \mathbb{R}_+^n} (J^{(5)} + J^{(6)} + J^{(7)}) \\ 
&\leq \lambda \sum_{2 \leq |\alpha| \leq d} \sum_{i,k=1}^n |h_{ik,\alpha}|^2 \, \varepsilon^{n-2} \, \int_{B_\delta(0) \cap \mathbb{R}_+^n} (\varepsilon + |x|)^{2|\alpha|+2-2n} \, dx \\ 
&+ C \, \delta^{2d+4-n} \, \varepsilon^{n-2}. 
\end{align*} 
Putting these facts together, the assertion follows. \\

\begin{proposition}
\label{technical.estimate.1}
If $\delta_0$ is sufficiently small, then we have 
\begin{align*}
&\int_{B_\delta(0) \cap \mathbb{R}_+^n} \Big ( u_\varepsilon^2 + \frac{n+2}{n-2} \, w^2 \Big )^{\frac{n}{n-2}} \\ 
&\leq \int_{B_\delta(0) \cap \mathbb{R}_+^n} (u_\varepsilon + w)^{\frac{2n}{n-2}} \\ 
&+ C \sum_{2 \leq |\alpha| \leq d} \sum_{i,k=1}^n |h_{ik,\alpha}|^2 \, \varepsilon^n \, \int_{B_\delta(0) \cap \mathbb{R}_+^n} (\varepsilon + |x|)^{2|\alpha|+2-2n} \, dx \\ 
&+ C \sum_{2 \leq |\alpha| \leq d} \sum_{i,k=1}^n |h_{ik,\alpha}| \, \delta^{|\alpha|-n} \, \varepsilon^n 
\end{align*} 
for all $0 < 2\varepsilon \leq \delta \leq \delta_0$.
\end{proposition}

\textbf{Proof.} 
The proof is analogous to the proof of Proposition 14 in \cite{Brendle2}. We omit the details. \\

\begin{proposition}
\label{technical.estimate.2} 
If $\delta_0$ is sufficiently small, then we have 
\begin{align*} 
&\int_{B_\delta(0) \cap \mathbb{R}_+^n} \frac{4(n-1)}{n-2} \, |du_\varepsilon|^2 + \int_{B_\delta(0) \cap \mathbb{R}_+^n} \frac{4(n-1)}{n-2} \, n(n+2) \, u_\varepsilon^{\frac{4}{n-2}} \, w^2 \\ 
&\leq Y(S_+^n,\partial S_+^n) \, \bigg ( \int_{B_\delta(0) \cap \mathbb{R}_+^n} (u_\varepsilon + w)^{\frac{2n}{n-2}} \bigg )^{\frac{n-2}{n}} \\ 
&+ \int_{\partial B_\delta(0) \cap \mathbb{R}_+^n} \frac{4(n-1)}{n-2} \sum_{i=1}^n \frac{x_i}{|x|} \, \partial_i u_\varepsilon \, u_\varepsilon \\ 
&+ C \sum_{2 \leq |\alpha| \leq d} \sum_{i,k=1}^n |h_{ik,\alpha}|^2 \, \varepsilon^n \, \int_{B_\delta(0) \cap \mathbb{R}_+^n} (\varepsilon + |x|)^{2|\alpha|+2-2n} \, dx \\ 
&+ C \sum_{2 \leq |\alpha| \leq d} \sum_{i,k=1}^n |h_{ik,\alpha}| \, \delta^{|\alpha|-n} \, \varepsilon^n 
\end{align*}
for all $0 < 2\varepsilon \leq \delta \leq \delta_0$.
\end{proposition}

\textbf{Proof.} 
The proof is similar to the proof of Proposition 15 in \cite{Brendle2}. We first observe that 
\[4n(n-1) \, \bigg ( \int_{\mathbb{R}_+^n} u_\varepsilon^{\frac{2n}{n-2}} \bigg )^{\frac{2}{n}} = Y(S_+^n,\partial S_+^n).\] 
Using H\"older's inequality, we obtain 
\begin{align*} 
&\int_{B_\delta(0) \cap \mathbb{R}_+^n} \frac{4(n-1)}{n-2} \, |du_\varepsilon|^2 - \int_{\partial B_\delta(0)} \frac{4(n-1)}{n-2} \sum_{i=1}^n \frac{x_i}{|x|} \, \partial_i u_\varepsilon \, u_\varepsilon \\ 
&+ \int_{B_\delta(0) \cap \mathbb{R}_+^n} \frac{4(n-1)}{n-2} \, n(n+2) \, u_\varepsilon^{\frac{4}{n-2}} \, w^2 \\ 
&= -\int_{B_\delta(0) \cap \mathbb{R}_+^n} \frac{4(n-1)}{n-2} \, \Delta u_\varepsilon \, u_\varepsilon + \int_{B_\delta(0) \cap \mathbb{R}_+^n} \frac{4(n-1)}{n-2} \, n(n+2) \, u_\varepsilon^{\frac{4}{n-2}} \, w^2 \\ 
&= \int_{B_\delta(0) \cap \mathbb{R}_+^n} 4n(n-1) \, u_\varepsilon^{\frac{4}{n-2}} \, \Big ( u_\varepsilon^2 + \frac{n+2}{n-2} \, w^2 \Big ) \\ 
&\leq Y(S_+^n,\partial S_+^n) \, \bigg ( \int_{B_\delta(0) \cap \mathbb{R}_+^n} \Big ( u_\varepsilon^2 + \frac{n+2}{n-2} \, w^2 \Big )^{\frac{n}{n-2}} \bigg )^{\frac{n-2}{n}}. 
\end{align*}
Hence, the assertion follows from Proposition \ref{technical.estimate.1}.

\section{Proof of the main result}

In this section, we construct a smooth function $v_{(\varepsilon,\delta)}: M \to \mathbb{R}$ with Yamabe energy less than $Y(S_+^n,\partial S_+^n)$. The existence of such a function is trivial when $Y(M,\partial M,g) \leq 0$. Hence, it suffices to consider the case $Y(M,\partial M,g) > 0$. As in the previous section, we fix a boundary point $p \in \partial M$. Moreover, we denote by $G: M \setminus \{p\} \to \mathbb{R}$ the Greens function for the conformal Laplacian with Neumann boundary condition with pole at $p$. In other words, $G$ satisfies 
\[\frac{4(n-1)}{n-2} \, \Delta_g G - R_g \, G = 0\] 
in $M \setminus \{p\}$ and $\partial_\nu G = 0$ along $\partial M \setminus \{p\}$. We assume that $G_p(x)$ is normalized so that $\lim_{x \to 0} |x|^{n-2} \, G(x) = 1$. With this normalization, we have 
\begin{equation} 
\label{estimate.for.Greens.function}
\big | G(x) - |x|^{2-n} \big | \leq C \sum_{2 \leq |\alpha| \leq d} \sum_{i,k=1}^n |h_{ik,\alpha}| \, |x|^{|\alpha|+2-n} + C \, |x|^{d+3-n}. 
\end{equation} 
Moreover, we consider the flux integral 
\begin{align*} 
\mathcal{I}(p,\delta) 
&= \int_{\partial B_\delta(0) \cap \mathbb{R}_+^n} \frac{4(n-1)}{n-2} \sum_{i=1}^n \frac{x_i}{|x|} \, (|x|^{2-n} \, \partial_i G - G \, \partial_i |x|^{2-n}) \\ 
&- \int_{\partial B_\delta(0) \cap \mathbb{R}_+^n} |x|^{1-2n} \sum_{i,k=1}^n x_i \, (|x|^2 \, \partial_k h_{ik} - 2n x_k \, h_{ik}) 
\end{align*}
where $\delta > 0$ is sufficiently small.

We next define a function $v_{(\varepsilon,\delta)}: M \to \mathbb{R}$ by 
\begin{equation} 
\label{definition.of.testfunction}
v_{(\varepsilon,\delta)} = \chi_\delta \, (u_\varepsilon + w) + (1 - \chi_\delta) \, \varepsilon^{\frac{n-2}{2}} \, G, 
\end{equation} 
where $\chi_\delta$ is the cut-off function defined above. Our main result is an upper bound for the Yamabe energy of $v_{(\varepsilon,\delta)}$:

\begin{proposition} 
\label{yamabe.energy.of.v} 
If $\delta_0$ is sufficiently small, then we have 
\begin{align*} 
&\int_M \Big ( \frac{4(n-1)}{n-2} \, |dv_{(\varepsilon,\delta)}|_g^2 + R_g \, v_{(\varepsilon,\delta)}^2 \Big ) \, d\text{\rm vol}_g \\ 
&\leq Y(S_+^n,\partial S_+^n) \, \bigg ( \int_M v_{(\varepsilon,\delta)}^{\frac{2n}{n-2}} \, d\text{\rm vol}_g \bigg )^{\frac{n-2}{n}} - \varepsilon^{n-2} \, \mathcal{I}(p,\delta) \\ 
&- \frac{\lambda}{2} \sum_{2 \leq |\alpha| \leq d} \sum_{i,k=1}^n |h_{ik,\alpha}|^2 \, \varepsilon^{n-2} \, \int_{B_\delta(0) \cap \mathbb{R}_+^n} (\varepsilon + |x|)^{2|\alpha|+2-2n} \, dx \\ 
&+ C \sum_{2 \leq |\alpha| \leq d} \sum_{i,k=1}^n |h_{ik,\alpha}| \, \delta^{|\alpha|+2-n} \, \varepsilon^{n-2} + C \, \delta^{2d+4-n} \, \varepsilon^{n-2} + C \, \delta^{-n} \, \varepsilon^n 
\end{align*} 
for all $0 < 2\varepsilon \leq \delta \leq \delta_0$.
\end{proposition} 

\textbf{Proof.} 
For abbreviation, we denote by $\Omega_\delta$ the set of all points in $M$ such that $x_1^2 + \hdots + x_n^2 < \delta^2$, where $(x_1,\hdots,x_n)$ denote the Fermi coordinates around $p$. (In other words, $\Omega_\delta$ is a coordinate ball, not a geodesic ball.) Using the divergence theorem, we obtain 
\begin{align*} 
&\int_{M \setminus \Omega_\delta} \Big ( \frac{4(n-1)}{n-2} \, |dv_{(\varepsilon,\delta)}|_g^2 + R_g \, v_{(\varepsilon,\delta)}^2 \Big ) \, d\text{\rm vol}_g \\ 
&= -\int_{M \setminus \Omega_\delta} \Big ( \frac{4(n-1)}{n-2} \, \Delta_g v_{(\varepsilon,\delta)} - R_g \, v_{(\varepsilon,\delta)} \Big ) \, \big ( v_{(\varepsilon,\delta)} - \varepsilon^{\frac{n-2}{2}} \, G \big ) \, d\text{\rm vol}_g \\ 
&- \int_{\partial \Omega_\delta} \frac{4(n-1)}{n-2} \, \partial_\nu v_{(\varepsilon,\delta)} \, v_{(\varepsilon,\delta)} \, d\sigma_g \\ 
&- \int_{\partial \Omega_\delta} \frac{4(n-1)}{n-2} \, \varepsilon^{\frac{n-2}{2}} \, (v_{(\varepsilon,\delta)} \, \partial_\nu G - G \, \partial_\nu v_{(\varepsilon,\delta)}) \, d\sigma_g, 
\end{align*} 
where $\nu$ denotes the outward-pointing unit normal to $\partial \Omega_\delta$. Note that 
\[v_{(\varepsilon,\delta)} - \varepsilon^{\frac{n-2}{2}} \, G = \chi_\delta \, \big ( u_\varepsilon + w - \varepsilon^{\frac{n-2}{2}} \, G \big )\] 
in $M \setminus \Omega_\delta$. In particular, we have $v_{(\varepsilon,\delta)} - \varepsilon^{\frac{n-2}{2}} \, G = 0$ in $M \setminus \Omega_{2\delta}$. Using (\ref{estimate.for.Greens.function}), we obtain  
\begin{align*} 
&\sup_{M \setminus \Omega_\delta} \big | v_{(\varepsilon,\delta)} - \varepsilon^{\frac{n-2}{2}} \, G \big | + \delta^2 \sup_{M \setminus \Omega_\delta} \Big | \frac{4(n-1)}{n-2} \, \Delta_g v_{(\varepsilon,\delta)} - R_g \, v_{(\varepsilon,\delta)} \Big | \\ 
&\leq C \sum_{2 \leq |\alpha| \leq d} \sum_{i,k=1}^n |h_{ik,\alpha}| \, \delta^{|\alpha|+2-n} \, \varepsilon^{\frac{n-2}{2}} + C \, \delta^{d+3-n} \, \varepsilon^{\frac{n-2}{2}} + C \, \delta^{-n} \, \varepsilon^{\frac{n+2}{2}}, 
\end{align*}
hence 
\begin{align*} 
&-\int_{M \setminus \Omega_\delta} \Big ( \frac{4(n-1)}{n-2} \, \Delta_g v_{(\varepsilon,\delta)} - R_g \, v_{(\varepsilon,\delta)} \Big ) \, \big ( v_{(\varepsilon,\delta)} - \varepsilon^{\frac{n-2}{2}} \, G \big ) \, d\text{\rm vol}_g \\ 
&\leq C \sum_{2 \leq |\alpha| \leq d} \sum_{i,k=1}^n |h_{ik,\alpha}|^2 \, \delta^{2|\alpha|+2-n} \, \varepsilon^{n-2} + C \, \delta^{2d+4-n} \, \varepsilon^{n-2} + C \, \delta^{-n-2} \, \varepsilon^{n+2}.  
\end{align*} 
We next observe that 
\begin{align*} 
&-\int_{\partial \Omega_\delta} \partial_\nu u_\varepsilon \, u_\varepsilon \, d\sigma_g \\ 
&\leq -\int_{\partial B_\delta(0) \cap \mathbb{R}_+^n} \sum_{i=1}^n \frac{x_i}{|x|} \, \partial_i u_\varepsilon \, u_\varepsilon + \int_{\partial B_\delta(0) \cap \mathbb{R}_+^n} \sum_{i=1}^n \frac{x_i}{|x|} \, u_\varepsilon \, \partial_k u_\varepsilon \, h_{ik} \\ 
&+ C \sum_{2 \leq |\alpha| \leq d} \sum_{i,k=1}^n |h_{ik,\alpha}|^2 \, \delta^{2|\alpha|+2-n} \, \varepsilon^{n-2} + C \, \delta^{2d+4-n} \, \varepsilon^{n-2}, 
\end{align*} 
hence 
\begin{align*} 
&-\int_{\partial \Omega_\delta} \partial_\nu v_{(\varepsilon,\delta)} \, v_{(\varepsilon,\delta)} \, d\sigma_g \\ 
&\leq -\int_{\partial B_\delta(0) \cap \mathbb{R}_+^n} \sum_{i=1}^n \frac{x_i}{|x|} \, \partial_i u_\varepsilon \, u_\varepsilon + \int_{\partial B_\delta(0) \cap \mathbb{R}_+^n} \sum_{i=1}^n \frac{x_i}{|x|} \, u_\varepsilon \, \partial_k u_\varepsilon \, h_{ik} \\ 
&+ C \sum_{2 \leq |\alpha| \leq d} \sum_{i,k=1}^n |h_{ik,\alpha}| \, \delta^{|\alpha|+2-n} \, \varepsilon^{n-2} + C \, \delta^{2d+4-n} \, \varepsilon^{n-2}. 
\end{align*} 
Moreover, we have 
\begin{align*} 
&-\int_{\partial \Omega_\delta} (v_{(\varepsilon,\delta)} \, \partial_\nu G - G \, \partial_\nu v_{(\varepsilon,\delta)}) \, d\sigma_g \\ 
&\leq -\int_{\partial B_\delta(0) \cap \mathbb{R}_+^n} \sum_{i=1}^n \frac{x_i}{|x|} \, (u_\varepsilon \, \partial_i G - G \, \partial_i u_\varepsilon) \\ 
&+ C \sum_{2 \leq |\alpha| \leq d} |h_{ik,\alpha}|^2 \, \delta^{2|\alpha|+2-n} \, \varepsilon^{\frac{n-2}{2}} + C \, \delta^{2d+4-n} \, \varepsilon^{\frac{n-2}{2}}. 
\end{align*} 
Putting these facts together, we obtain 
\begin{align*} 
&\int_{M \setminus \Omega_\delta} \Big ( \frac{4(n-1)}{n-2} \, |dv_{(\varepsilon,\delta)}|_g^2 + R_g \, v_{(\varepsilon,\delta)}^2 \Big ) \, d\text{\rm vol}_g \\ 
&\leq -\int_{\partial \Omega_\delta} \frac{4(n-1)}{n-2} \sum_{i=1}^n \frac{x_i}{|x|} \, \partial_i u_\varepsilon \, u_\varepsilon + \int_{\partial B_\delta(0) \cap \mathbb{R}_+^n} \frac{4(n-1)}{n-2} \sum_{i=1}^n \frac{x_i}{|x|} \, u_\varepsilon \, \partial_k u_\varepsilon \, h_{ik} \\ 
&- \int_{\partial B_\delta(0) \cap \mathbb{R}_+^n} \frac{4(n-1)}{n-2} \, \varepsilon^{\frac{n-2}{2}} \sum_{i=1}^n \frac{x_i}{|x|} \, (u_\varepsilon \, \partial_i G - G \, \partial_i u_\varepsilon) \\ 
&+ C \sum_{2 \leq |\alpha| \leq d} \sum_{i,k=1}^n |h_{ik,\alpha}| \, \delta^{|\alpha|+2-n} \, \varepsilon^{n-2} + C \, \delta^{2d+4-n} \, \varepsilon^{n-2} + C \, \delta^{-n-2} \, \varepsilon^{n+2}. 
\end{align*} 
On the other hand, it follows from Proposition \ref{estimate.for.test.function} and Proposition \ref{technical.estimate.2} that 
\begin{align*} 
&\int_{\Omega_\delta} \Big ( \frac{4(n-1)}{n-2} \, |dv_{(\varepsilon,\delta)}|_g^2 + R_g \, v_{(\varepsilon,\delta)}^2 \Big ) \, d\text{\rm vol}_g \\ 
&\leq Y(S_+^n,\partial S_+^n) \, \bigg ( \int_{\Omega_\delta} v_{(\varepsilon,\delta)}^{\frac{2n}{n-2}} \, d\text{\rm vol}_g \bigg )^{\frac{n-2}{n}} + \int_{\partial B_\delta(0) \cap \mathbb{R}_+^n} \frac{4(n-1)}{n-2} \sum_{i=1}^n \frac{x_i}{|x|} \, \partial_i u_\varepsilon \, u_\varepsilon \\ 
&+ \int_{\partial B_\delta(0) \cap \mathbb{R}_+^n} \sum_{i,k=1}^n \frac{x_i}{|x|} \, (u_\varepsilon^2 \, \partial_k h_{ik} - 2 \, u_\varepsilon \, \partial_k u_\varepsilon \, h_{ik}) \\ 
&- \frac{\lambda}{2} \sum_{2 \leq |\alpha| \leq d} \sum_{i,k=1}^n |h_{ik,\alpha}|^2 \, \varepsilon^{n-2} \, \int_{B_\delta(0) \cap \mathbb{R}_+^n} (\varepsilon + |x|)^{2|\alpha|+2-2n} \, dx \\ 
&+ C \sum_{2 \leq |\alpha| \leq d}  \sum_{i,k=1}^n |h_{ik,\alpha}| \, \delta^{|\alpha|+2-n} \, \varepsilon^{n-2} + C \, \delta^{2d+4-n} \, \varepsilon^{n-2}. 
\end{align*} 
If we add the last two inequalities, we obtain 
\begin{align*} 
&\int_M \Big ( \frac{4(n-1)}{n-2} \, |dv_{(\varepsilon,\delta)}|_g^2 + R_g \, v_{(\varepsilon,\delta)}^2 \Big ) \, d\text{\rm vol}_g \\ 
&\leq Y(S_+^n,\partial S_+^n) \, \bigg ( \int_{\Omega_\delta} v_{(\varepsilon,\delta)}^{\frac{2n}{n-2}} \, d\text{\rm vol}_g \bigg )^{\frac{n-2}{n}} \\ 
&+ \int_{\partial B_\delta(0) \cap \mathbb{R}_+^n} \sum_{i,k=1}^n \frac{x_i}{|x|} \, \Big ( u_\varepsilon^2 \, \partial_k h_{ik} + \frac{2n}{n-2} \, u_\varepsilon \, \partial_k u_\varepsilon \, h_{ik} \Big ) \\ 
&- \int_{\partial B_\delta(0) \cap \mathbb{R}_+^n} \frac{4(n-1)}{n-2} \, \varepsilon^{\frac{n-2}{2}} \sum_{i=1}^n \frac{x_i}{|x|} \, (u_\varepsilon \, \partial_i G - G \, \partial_i u_\varepsilon) \\ 
&- \frac{\lambda}{2} \sum_{2 \leq |\alpha| \leq d} \sum_{i,k=1}^n |h_{ik,\alpha}|^2 \, \varepsilon^{n-2} \, \int_{B_\delta(0) \cap \mathbb{R}_+^n} (\varepsilon + |x|)^{2|\alpha|+2-2n} \, dx \\ 
&+ C \sum_{2 \leq |\alpha| \leq d}  \sum_{i,k=1}^n |h_{ik,\alpha}| \, \delta^{|\alpha|+2-n} \, \varepsilon^{n-2} + C \, \delta^{2d+4-n} \, \varepsilon^{n-2} + C \, \delta^{-n-2} \, \varepsilon^{n+2}. 
\end{align*} From this, the assertion follows easily. \\

\begin{theorem}
Assume that $p \notin \mathcal{Z}$. Then $Y(M,\partial M,g) < Y(S_+^n,\partial S_+^n)$.
\end{theorem}

\textbf{Proof.} 
Since $p \notin \mathcal{Z}$, we have $\sum_{2 \leq |\alpha| \leq d} \sum_{i,k=1}^n |h_{ik,\alpha}|^2 > 0$. Using Proposition \ref{yamabe.energy.of.v}, we obtain 
\begin{align*} 
&\int_M \Big ( \frac{4(n-1)}{n-2} \, |dv_{(\varepsilon,\delta)}|_g^2 + R_g \, v_{(\varepsilon,\delta)}^2 \Big ) \, d\text{\rm vol}_g \\ 
&< Y(S_+^n,\partial S_+^n) \, \bigg ( \int_M v_{(\varepsilon,\delta)}^{\frac{2n}{n-2}} \, d\text{\rm vol}_g \bigg )^{\frac{n-2}{n}} 
\end{align*}
if $\varepsilon > 0$ is sufficiently small. From this, the assertion follows. \\

In the remainder of this section, we study the case $p \in \mathcal{Z}$. In this case, we consider the manifold $(M \setminus \{p\},G^{\frac{4}{n-2}} \, g)$. This manifold is scalar flat and its boundary is totally geodesic. After doubling this manifold, we obtain an asymptotically flat manifold with zero scalar curvature.

\begin{proposition}
\label{mass}
Assume that $p \in \mathcal{Z}$. Then the following statements hold: 

(i) The limit $\lim_{\delta \to 0} \mathcal{I}(p,\delta)$ exists.

(ii) The doubling of $(M \setminus \{p\},G^{\frac{4}{n-2}} \, g)$ has a well-defined mass which equals $\lim_{\delta \to 0} \mathcal{I}(p,\delta)$ up to a positive factor.
\end{proposition}

\textbf{Proof.} 
For abbreviation, let $\overline{g} = G^{\frac{4}{n-2}} \, g$. We consider the inverted coordinates $y = \frac{x}{|x|^2}$, where $(x_1,\hdots,x_n)$ are conformal Fermi coordinates around $p$. In these coordinates, the metric $\overline{g}$ is given by 
\begin{align*} 
\overline{g} \Big ( \frac{\partial}{\partial y_j},\frac{\partial}{\partial y_l} \Big ) 
&= \Big [ 1 + \Phi \Big ( \frac{y}{|y|^2} \Big ) \Big ]^{\frac{4}{n-2}} \\ 
&\cdot |y|^{-4} \sum_{i,k=1}^n (|y|^2 \, \delta_{ij} - 2y_iy_j) \, (|y|^2 \, \delta_{kl} - 2y_ky_l) \, g_{ik} \Big ( \frac{y}{|y|^2} \Big ), 
\end{align*}
where $\Phi(x) = |x|^{n-2} \, G(x) - 1$. Using the relations $g_{ik}(x) = \delta_{ik} + h_{ik}(x) + O(|x|^{2d+2})$ and $\Phi(x) = O(|x|^{d+1})$, we obtain 
\begin{align*} 
\overline{g} \Big ( \frac{\partial}{\partial y_j},\frac{\partial}{\partial y_l} \Big ) 
&= \Big [ 1 + \frac{4}{n-2} \, \Phi \Big ( \frac{y}{|y|^2} \Big ) \Big ] \, \delta_{jl} \\ 
&+ |y|^{-4} \sum_{i,k=1}^n (|y|^2 \, \delta_{ij} - 2y_iy_j) \, (|y|^2 \, \delta_{kl} - 2y_ky_l) \, h_{ik} \Big ( \frac{y}{|y|^2} \Big ) \\ 
&+ O(|y|^{-2d-2}). 
\end{align*}
In particular, we have $\overline{g} \big ( \frac{\partial}{\partial y_j},\frac{\partial}{\partial y_l} \big ) = \delta_{jl} + O(|y|^{-d-1})$. Hence, the doubling of $(M \setminus \{p\},\overline{g})$ is asymptotically flat in the sense of Bartnik \cite{Bartnik}, and has a well-defined ADM mass.

Since $\text{\rm tr} \, h = O(|x|^{2d+2})$, it follows that 
\[\sum_{j,l=1}^n y_j \, \frac{\partial}{\partial y_j} \overline{g} \Big ( \frac{\partial}{\partial y_l},\frac{\partial}{\partial y_l} \Big ) = -\frac{4n}{n-2} \, \sum_{j=1}^n \frac{y_j}{|y|^2} \, (\partial_j \Phi) \Big ( \frac{y}{|y|^2} \Big ) + O(|y|^{-2d-2}).\] 
Moreover, we have 
\begin{align*} 
\sum_{j,l=1}^n y_j \, \frac{\partial}{\partial y_l} \overline{g} \Big ( \frac{\partial}{\partial y_j},\frac{\partial}{\partial y_l} \Big ) 
&= -\frac{4}{n-2} \sum_{i=1}^n \frac{y_i}{|y|^2} \, (\partial_i \Phi) \Big ( \frac{y}{|y|^2} \Big ) \\ 
&- \sum_{i,k=1}^n \frac{y_i}{|y|^2} \, (\partial_k h_{ik}) \Big ( \frac{y}{|y|^2} \Big ) \\ 
&+ 2n \sum_{i,k=1}^n \frac{y_i \, y_k}{|y|^2} \, h_{ik} \Big ( \frac{y}{|y|^2} \Big ) + O(|y|^{-2d-2}), 
\end{align*}
where $\partial_i \Phi(x) = \frac{\partial}{\partial x_i} \Phi(x)$. Putting these facts together, we obtain 
\begin{align*} 
&\sum_{j,l=1}^n y_j \, \frac{\partial}{\partial y_l} \overline{g} \Big ( \frac{\partial}{\partial y_j},\frac{\partial}{\partial y_l} \Big ) - \sum_{j,l=1}^n y_j \, \frac{\partial}{\partial y_j} \overline{g} \Big ( \frac{\partial}{\partial y_l},\frac{\partial}{\partial y_l} \Big ) \\ 
&= \frac{4(n-1)}{n-2} \sum_{i=1}^n \frac{y_i}{|y|^2} \, (\partial_i \Phi) \Big ( \frac{y}{|y|^2} \Big ) \\ 
&- \sum_{i,k=1}^n \frac{y_i}{|y|^2} \, (\partial_k h_{ik}) \Big ( \frac{y}{|y|^2} \Big ) + 2n \sum_{i,k=1}^n \frac{y_i \, y_k}{|y|^2} \, h_{ik} \Big ( \frac{y}{|y|^2} \Big ) \\ 
&+ O(|y|^{-2d-2}). 
\end{align*}
This implies 
\begin{align*} 
&\int_{\partial B_{\frac{1}{\delta}}(0) \cap \mathbb{R}_+^n} \sum_{j,l=1}^n \frac{y_j}{|y|} \, \frac{\partial}{\partial y_l} \overline{g} \Big ( \frac{\partial}{\partial y_j},\frac{\partial}{\partial y_l} \Big ) - \int_{\partial B_{\frac{1}{\delta}}(0) \cap \mathbb{R}_+^n} \sum_{j,l=1}^n \frac{y_j}{|y|} \, \frac{\partial}{\partial y_j} \overline{g} \Big ( \frac{\partial}{\partial y_l},\frac{\partial}{\partial y_l} \Big ) \\ 
&= \int_{\partial B_\delta(0) \cap \mathbb{R}_+^n} \frac{4(n-1)}{n-2} \, |x|^{3-2n} \sum_{i=1}^n x_i \, \partial_i \Phi(x) \\ 
&- \int_{\partial B_\delta(0) \cap \mathbb{R}_+^n} |x|^{3-2n} \sum_{i,k=1}^n x_i \, (\partial_k h_{ik})(x) + \int_{\partial B_\delta(0) \cap \mathbb{R}_+^n} 2n \, |x|^{1-2n} \sum_{i,k=1}^n x_i \, x_k \, h_{ik}(x) \\ 
&+ O(\delta^{2d+4-n}) \\ 
&= \mathcal{I}(p,\delta) + O(\delta^{2d+n-4}).
\end{align*}
As $\delta \to 0$, the left hand side converges to a positive multiple of the ADM mass. From this the assertion follows. \\

\begin{theorem}
\label{degenerate.case}
Assume that $p \in \mathcal{Z}$. If $\lim_{\delta \to 0} \mathcal{I}(p,\delta)$ is positive, then $Y(M,\partial M,g) < Y(S_+^n,\partial S_+^n)$.
\end{theorem}

\textbf{Proof.}
Since $p \in \mathcal{Z}$, we have $\sum_{2 \leq |\alpha| \leq d} \sum_{i,k=1}^n |h_{ik,\alpha}|^2 = 0$. By Proposition \ref{yamabe.energy.of.v}, we can find positive real numbers $\delta_0$ and $C$ such that 
\begin{align*} 
&\int_M \Big ( \frac{4(n-1)}{n-2} \, |dv_{(\varepsilon,\delta)}|_g^2 + R_g \, v_{(\varepsilon,\delta)}^2 \Big ) \, d\text{\rm vol}_g \\ 
&\leq Y(S_+^n,\partial S_+^n) \, \bigg ( \int_M v_{(\varepsilon,\delta)}^{\frac{2n}{n-2}} \, d\text{\rm vol}_g \bigg )^{\frac{n-2}{n}} - \varepsilon^{n-2} \, \mathcal{I}(p,\delta) \\ 
&+ C \, \delta^{2d+4-n} \, \varepsilon^{n-2} + C \, \delta^{-n} \, \varepsilon^n 
\end{align*} 
whenever $0 < 2\varepsilon \leq \delta \leq \delta_0$. Since $\lim_{\delta \to 0} \mathcal{I}(p,\delta)$ is positive, we can find a real number $\delta \in (0,\delta_0]$ such that $\mathcal{I}(p,\delta) > C \, \delta^{2d+4-n}$. In the next step, we choose $\varepsilon \in (0,\frac{\delta}{2}]$ small enough so that $\mathcal{I}(p,\delta) > C \, \delta^{2d+4-n} + C \, \delta^{-n} \, \varepsilon^2$. For this choice of $\varepsilon$ and $\delta$, we have 
\begin{align*} 
&\int_M \Big ( \frac{4(n-1)}{n-2} \, |dv_{(\varepsilon,\delta)}|_g^2 + R_g \, v_{(\varepsilon,\delta)}^2 \Big ) \, d\text{\rm vol}_g \\ 
&< Y(S_+^n,\partial S_+^n) \, \bigg ( \int_M v_{(\varepsilon,\delta)}^{\frac{2n}{n-2}} \, d\text{\rm vol}_g \bigg )^{\frac{n-2}{n}}. 
\end{align*}
This completes the proof.

\appendix

\section{An elliptic system on $\mathbb{R}_+^n$}

In this section, we describe the construction of the vector field $V$. In the following, we consider the hemisphere $S_+^n$, equipped with the round metric of constant sectional curvature $4$. We denote by $\mathcal{X}$ the space of all vector fields $V$ on $S_+^n$ such that $V$ is of class $H^1$ and $\langle V,\nu \rangle = 0$ along $\partial S_+^n$. Moreover, we denote by $\mathcal{Y}$ the space of all trace-free symmetric two-tensors on $S_+^n$ of class $L^2$. We next define a linear operator $\mathcal{D}: \mathcal{X} \to \mathcal{Y}$ by
\[\mathcal{D} V = \widehat{\mathscr{L}}_V g = \mathscr{L}_V g - \frac{2}{n} \, (\text{\rm div}_g \, V) \, g.\]
In other words, $\mathcal{D}$ is the conformal Killing operator.

\begin{lemma}
\label{estim.1}
We have
\[\|\nabla V\|_{L^2(S_+^n)}^2 \leq \|\mathcal{D} V\|_{L^2(S_+^n)}^2 + 4(n-1) \, \|V\|_{L^2(S_+^n)}^2\]
for all $V \in \mathcal{X}$.
\end{lemma}

\textbf{Proof.}
Without loss of generality, we may assume that $V$ is smooth. By definition of $\mathcal{D}$, we have
\[\|\mathcal{D} V\|_{L^2(S_+^n)}^2 = \int_{S_+^n} \Big [ \nabla_i V^k \, \nabla^i V_k + \nabla_i V^k \, \nabla_k V^i - \frac{2}{n} \, (\text{\rm div}_g \, V)^2 \Big ] \, d\text{\rm vol}_g.\]
Integration by parts yields 
\begin{align*}
\int_{S_+^n} \nabla_i V^k \, \nabla_k V^i \, d\text{\rm vol}_g
&= -\int_{S_+^n} V^k \, \nabla_i \nabla_k V^i \, d\text{\rm vol}_g \\ 
&= -\int_{S_+^n} V^k \, \nabla_k \nabla_i V^i \, d\text{\rm vol}_g - \int_{S_+^n} \text{\rm Ric}_{ik} \, V^i \, V^k \, d\text{\rm vol}_g \\ 
&= \int_{S_+^n} (\text{\rm div}_g \, V)^2 \, d\text{\rm vol}_g - 4(n-1) \int_{S_+^n} |V|^2 \, d\text{\rm vol}_g.
\end{align*}
Putting these facts together, we obtain
\[\|\mathcal{D} V\|_{L^2(S_+^n)}^2 + 4(n-1) \, \|V\|_{L^2(S_+^n)}^2 = \|\nabla V\|_{L^2(S_+^n)}^2 + \frac{n-2}{n} \, \|\text{\rm div}_g \, V\|_{L^2(S_+^n)}^2.\] From this, the assertion follows. \\

It follows from Lemma \ref{estim.1} and Rellich's theorem that $\ker \mathcal{D}$ is finite-dimensional. We now consider the subspace
\[\mathcal{X}_0 = \{V \in \mathcal{X}: \text{\rm $\langle V,W \rangle_{L^2(S_+^n)} = 0$ for all $W \in \ker \mathcal{D}$}\}.\] 

\begin{lemma}
\label{estim.2}
We have
\[\|V\|_{L^2(S_+^n)}^2 + \|\nabla V\|_{L^2(S_+^n)}^2 \leq K \, \|\mathcal{D} V\|_{L^2(S_+^n)}^2\]
for all $V \in \mathcal{X}_0$. Here, $K$ is a positive constant that depends only on $n$.
\end{lemma}

\textbf{Proof.} 
Suppose that the assertion is false. Then we can find a sequence of vector fields $V^{(\nu)} \in \mathcal{X}_0$ such that 
\begin{equation} 
\label{contradiction}
\|V^{(\nu)}\|_{L^2(S_+^n)}^2 + \|\nabla V^{(\nu)}\|_{L^2(S_+^n)}^2 = 1 
\end{equation}
for all $\nu$ and $\|\mathcal{D} V^{(\nu)}\|_{L^2(S_+^n)} \to 0$ as $\nu \to \infty$. After passing to a subsequence, we may assume that the sequence $V^{(\nu)}$ converges weakly to a vector field $W \in \mathcal{X}_0$. Then $\mathcal{D} W = 0$. Since $W \in \mathcal{X}_0$, we conclude that $W = 0$. This implies $\|V^{(\nu)}\|_{L^2(S_+^n)} \to 0$ as $\nu \to \infty$. Using Lemma \ref{estim.1}, we obtain $\|\nabla V^{(\nu)}\|_{L^2(S_+^n)} \to 0$ as $\nu \to \infty$. This contradicts (\ref{contradiction}). \\

\begin{proposition}
\label{existence.of.weak.solution}
Given any $h \in \mathcal{Y}$, there exists a unique vector field $V \in \mathcal{X}_0$ such that $\langle h - \mathcal{D} V,\mathcal{D} W \rangle_{L^2(S_+^n)} = 0$ for all $W \in \mathcal{X}$. The vector field $V$ satisfies the estimate 
\begin{equation} 
\label{estimate.for.V.1}
\|V\|_{L^2(S_+^n)}^2 + \|\nabla V\|_{L^2(S_+^n)}^2 \leq K \, \|h\|_{L^2(S_+^n)}^2. 
\end{equation}
\end{proposition}

\textbf{Proof.} 
It follows from Lemma \ref{estim.2} that the operator $\mathcal{D}: \mathcal{X}_0 \to \mathcal{Y}$ has closed range. Hence, we can find a vector field $V \in \mathcal{X}_0$ such that $\|h - \mathcal{D} V\|_{L^2(S_+^n)}^2$ is minimal. This vector field satisfies $\langle h - \mathcal{D} V,\mathcal{D} W \rangle_{L^2(S_+^n)} = 0$ for all $W \in \mathcal{X}_0$. This proves the existence statement.

We next assume that $V \in \mathcal{X}_0$ is a vector field satisfying $\langle h - \mathcal{D} V,\mathcal{D} W \rangle_{L^2(S_+^n)} = 0$ for all $W \in \mathcal{X}_0$. This implies $\langle h - \mathcal{D} V,\mathcal{D} V \rangle = 0$, hence $\|\mathcal{D} V\|_{L^2(S_+^n)}^2 \leq \|h\|_{L^2(S_+^n)}^2$. Thus, we conclude that 
\[\|V\|_{L^2(S_+^n)}^2 + \|\nabla V\|_{L^2(S_+^n)}^2 \leq K \, \|\mathcal{D} V\|_{L^2(S_+^n)}^2 \leq K \, \|h\|_{L^2(S_+^n)}^2\] 
by Lemma \ref{estim.2}. In particular, if $h = 0$, then $V = 0$. From this, the uniqueness statement follows. \\

In the next step, we consider the stereographic projection from $S_+^n$ to $\mathbb{R}_+^n \cup \{\infty\}$. The metric $g$ can be written in the form $g_{ik} = u^{\frac{4}{n-2}} \, \delta_{ik}$, where 
\[u(x) = \Big ( \frac{1}{1 + |x|^2} \Big )^{\frac{n-2}{2}}.\]

\begin{theorem}
\label{existence.of.smooth.solution}
Let $h$ be a trace-free symmetric two-tensor on $\mathbb{R}_+^n$. We assume that $h$ is smooth and has compact support. Then there exists a smooth vector field $V$ on $\mathbb{R}_+^n$ with the following properties: 
\begin{itemize}
\item At each point $x \in \mathbb{R}_+^n$, we have 
\[\sum_{k=1}^n \partial_k \Big [ u^{\frac{2n}{n-2}} \, \big ( h_{ik} - \partial_i V_k - \partial_k V_i + \frac{2}{n} \, \text{\rm div} \, V \, \delta_{ik} \big ) \Big ] = 0\] 
for all $i \in \{1,\hdots,n\}$. 
\item At each point $x \in \partial \mathbb{R}_+^n$, we have $V_n(x) = \partial_n V_i(x) - h_{in}(x) = 0$ for all $i \in \{1,\hdots,n-1\}$.
\end{itemize}
Moreover, the vector field $V$ satisfies  
\begin{equation}
\label{estimate.for.V.2}
\int_{\mathbb{R}_+^n} u(x)^{\frac{2(n+2)}{n-2}} \, |V(x)|^2 \, dx \leq K \int_{\mathbb{R}_+^n} u(x)^{\frac{2n}{n-2}} \, |h(x)|^2 \, dx. 
\end{equation} 
\end{theorem}

\textbf{Proof.} 
By Proposition \ref{existence.of.weak.solution}, there exists a smooth vector field $V \in \mathcal{X}_0$ such that 
\[\int_{\mathbb{R}_+^n} \langle u^{\frac{4}{n-2}} \, h - \mathcal{D} V,\mathcal{D} W \rangle_g \, d\text{\rm vol}_g = 0\] 
for all vector fields $W \in \mathcal{X}$. This implies
\begin{equation} 
\label{weak.formulation}
\int_{\mathbb{R}_+^n} u^{\frac{2n}{n-2}} \sum_{i,k=1}^n \big ( h_{ik} - \partial_i V_k - \partial_k V_i + \frac{2}{n} \, \text{\rm div} \, V \, \delta_{ik} \big ) \, \partial_k W_i \, dx = 0 
\end{equation}
for all $W \in \mathcal{X}$. Since $V \in \mathcal{X}_0$, we have $V_n(x) = 0$ for $x \in \partial \mathbb{R}_+^n$.

By assumption, $h$ is smooth. Using general regularity results for elliptic systems (cf. \cite{Hormander}, \cite{Taylor}), we conclude that $V$ is smooth. Using (\ref{weak.formulation}), we obtain
\[\sum_{k=1}^n \partial_k \Big [ u^{\frac{2n}{n-2}} \, \big ( h_{ik} - \partial_i V_k - \partial_k V_i + \frac{2}{n} \, \text{\rm div} \, V \, \delta_{ik} \big ) \Big ] = 0\] 
for all points $x \in \mathbb{R}_+^n$ and all $i \in \{1,\hdots,n\}$. Moreover, at each point $x \in \partial \mathbb{R}_+^n$, we have $\partial_n V_i(x) - h_{in}(x) = 0$ for $i \in \{1,\hdots,n-1\}$. Finally, the estimate (\ref{estimate.for.V.2}) follows immediately from (\ref{estimate.for.V.1}). \\

\begin{proposition} 
\label{weighted.estimate}
Fix a real number $\sigma$ such that $1 < \sigma < n-2$. Let $h$ be a trace-free symmetric two-tensor on $\mathbb{R}_+^n$ which is smooth and has compact support. Moreover, let $V$ be the vector field constructed in Theorem \ref{existence.of.smooth.solution}. Finally, let us assume that 
\[\sup_{r \geq 1} r^{-2\sigma-n-2} \int_{(B_{2r}(0) \setminus B_r(0)) \cap \mathbb{R}_+^n} |V(x)|^2 \, dx < \infty.\] 
Then there exists a constant $C$, depending only on $n$ and $\sigma$, such that 
\begin{align} 
\label{weighted.estimate.for.V}
&\sup_{r \geq 1} r^{-2\sigma-n-2} \int_{(B_{2r}(0) \setminus B_r(0)) \cap \mathbb{R}_+^n} |V(x)|^2 \, dx \notag \\ 
&\leq C \int_{\mathbb{R}_+^n} (1 + |x|^2)^{-n-2} \, |V(x)|^2 \, dx \\ 
&+ C \sup_{r \geq 1} r^{-2\sigma-n} \int_{(B_{2r}(0) \setminus B_r(0)) \cap \mathbb{R}_+^n} |h(x)|^2 \, dx. \notag 
\end{align} 
\end{proposition}

\textbf{Proof.}
We extend $V$ and $h$ to $\mathbb{R}^n$ by reflection. More precisely, we define a vector field $\tilde{V}$ on $\mathbb{R}^n$ by 
\begin{align*}
&\tilde{V}_i(x_1,\hdots,x_{n-1},x_n) = \tilde{V}_i(x_1,\hdots,x_{n-1},-x_n) = V_i(x_1,\hdots,x_{n-1},x_n) \\
&\tilde{V}_n(x_1,\hdots,x_{n-1},x_n) = -\tilde{V}_n(x_1,\hdots,x_{n-1},-x_n) = V_n(x_1,\hdots,x_{n-1},x_n)
\end{align*}
for all $x \in \mathbb{R}_+^n$ and all $i \in \{1,\hdots,n-1\}$. Similarly, we define a trace-free symmetric two-tensor $\tilde{h}$ on $\mathbb{R}^n$ by 
\begin{align*}
&\tilde{h}_{ik}(x_1,\hdots,x_{n-1},x_n) = \tilde{h}_{ik}(x_1,\hdots,x_{n-1},-x_n) = h_{ik}(x_1,\hdots,x_{n-1},x_n) \\
&\tilde{h}_{in}(x_1,\hdots,x_{n-1},x_n) = -\tilde{h}_{in}(x_1,\hdots,x_{n-1},-x_n) = h_{in}(x_1,\hdots,x_{n-1},x_n) \\
&\tilde{h}_{nk}(x_1,\hdots,x_{n-1},x_n) = -\tilde{h}_{nk}(x_1,\hdots,x_{n-1},-x_n) = h_{nk}(x_1,\hdots,x_{n-1},x_n) \\
&\tilde{h}_{nn}(x_1,\hdots,x_{n-1},x_n) = \tilde{h}_{nn}(x_1,\hdots,x_{n-1},-x_n) = h_{nn}(x_1,\hdots,x_{n-1},x_n)
\end{align*}
for all $x \in \mathbb{R}_+^n$ and all $i,k \in \{1,\hdots,n-1\}$.

Since $V \in \mathcal{X}$, we have $V_n(x) = 0$ for all $x \in \partial \mathbb{R}_+^n$. Consequently, the vector field $\tilde{V}$ is a vector field on $S^n$ of class $H^1$. We claim that
\begin{equation}
\label{weak.formulation.V.tilde}
\int_{\mathbb{R}^n} u^{\frac{2n}{n-2}} \sum_{i,k=1} \big ( \tilde{h}_{ik} - \partial_i \tilde{V}_k - \partial_k \tilde{V}_i + \frac{2}{n} \, \text{\rm div} \, \tilde{V} \, \delta_{ik} \big ) \, \partial_k \tilde{W}_i \, dx = 0
\end{equation}
for all vector fields $\tilde{W}$ on $S^n$ of class $H^1$. In order to prove (\ref{weak.formulation.V.tilde}), we fix a vector field $\tilde{W}$ of class $H^1$. We then define a vector field $W$ on $S_+^n$ by
\begin{align*}
&W_i(x_1,\hdots,x_{n-1},x_n) = \tilde{W}_i(x_1,\hdots,x_{n-1},x_n) + \tilde{W}_i(x_1,\hdots,x_{n-1},-x_n) \\ 
&W_n(x_1,\hdots,x_{n-1},x_n) = \tilde{W}_n(x_1,\hdots,x_{n-1},x_n) - \tilde{W}_n(x_1,\hdots,x_{n-1},-x_n)
\end{align*} 
for all $x \in \mathbb{R}_+^n$ and all $i \in \{1,\hdots,n-1\}$. Clearly, $W \in \mathcal{X}$. Therefore, we have
\[\int_{\mathbb{R}_+^n} u^{\frac{2n}{n-2}} \sum_{i,k=1}^n \big ( h_{ik} - \partial_i V_k - \partial_k V_i + \frac{2}{n} \, \text{\rm div} \, V \, \delta_{ik} \big ) \, \partial_k W_i \, dx = 0\] 
by definition of $V$. From this, the identity (\ref{weak.formulation.V.tilde}) follows easily.

We now complete the proof of Proposition \ref{weighted.estimate}. Using Proposition 23 in \cite{Brendle2}, we obtain
\begin{align*} 
&\sup_{r \geq 1} r^{-2\sigma-n-2} \int_{B_{2r}(0) \setminus B_r(0)} |\tilde{V}(x)|^2 \, dx \\ 
&\leq C \int_{\mathbb{R}_+^n} (1 + |x|^2)^{-n-2} \, |\tilde{V}(x)|^2 \, dx \\ 
&+ C \sup_{r \geq 1} r^{-2\sigma-n} \int_{B_{2r}(0) \setminus B_r(0)} |\tilde{h}(x)|^2 \, dx. 
\end{align*} 
Here, $C$ is a positive constant that depends only on $\sigma$ and $n$. (In \cite{Brendle2}, this result was stated in the special case that $\tilde{V}$ and $\tilde{h}$ are smooth, but the proof only requires that $\tilde{h}$ belongs to $L^2$ and $\tilde{V}$ is of class $H^1$.) From this the assertion follows. \\

\begin{corollary} 
\label{polynomial}
Consider a trace-free symmetric two-tensor of the form 
\[h_{ik}(x) = \chi(|x|/\rho) \, \sum_{2 \leq |\alpha| \leq d} h_{ik,\alpha} \, x^\alpha,\] 
where $d = [\frac{n-2}{2}]$, $\rho \geq 1$, and $\chi: \mathbb{R} \to \mathbb{R}$ is a fixed cutoff function satisfying $\chi(t) = 0$ for $t \geq 2$. Let $V$ be the vector field constructed in Theorem \ref{existence.of.smooth.solution}. Then, for every multi-index $\beta$, we have 
\begin{equation} 
|\partial^\beta V(x)|^2 \leq C \sum_{2 \leq |\alpha| \leq d} |h_{ik,\alpha}|^2 \, (1 + |x|^2)^{|\alpha|+1-|\beta|} 
\end{equation}
for all $x \in \mathbb{R}_+^n$. Here, $C$ is positive constant which depends on $n$ and $|\beta|$, but not on $\rho$.
\end{corollary}

\textbf{Proof.} Without loss of generality, we may assume that 
\[h_{ik}(x) = \chi(|x|/\rho) \, \sum_{|\alpha| = d'} h_{ik,\alpha} \, x^\alpha,\] 
where $2 \leq d' \leq d$. Since $d' < \frac{n}{2}$, we have 
\[\int_{\mathbb{R}_+^n} (1 + |x|^2)^{-n} \, |h(x)|^2 \, dx \leq C \sum_{|\alpha| = d'} |h_{ik,\alpha}|^2\] 
for some uniform constant $C$. Using (\ref{estimate.for.V.2}), we obtain  
\[\int_{\mathbb{R}_+^n} (1 + |x|^2)^{-n-2} \, |V(x)|^2 \, dx \leq C \sum_{|\alpha| = d'} |h_{ik,\alpha}|^2.\] 
We now apply Proposition \ref{weighted.estimate} with $\sigma = d'$. This yields 
\[\sup_{r \geq 1} r^{-2d'-n-2} \int_{\{r \leq |x| \leq 2r\}} |V(x)|^2 \, dx \leq C \sum_{|\alpha| = d'} |h_{ik,\alpha}|^2.\] 
Using elliptic estimates, we conclude that  
\[|\partial^\beta V(x)|^2 \leq C \sum_{|\alpha| = d'} |h_{ik,\alpha}|^2 \, (1 + |x|^2)^{d'+1-|\beta|}\] 
for every multi-index $\beta$.

\end{document}